\documentclass{m2an}
\usepackage[colorlinks]{hyperref}
\usepackage{amsmath}
\usepackage{amsthm}
\usepackage{amssymb}
\usepackage[utf8]{inputenc}
\usepackage{chngcntr}
\usepackage{apptools}
\usepackage{graphicx}

\newcommand{\dd}{{\rm{d}}}
\def\supp{\textup{supp}}
\newcommand{\ii}{{\rm{i}}}
\usepackage{xcolor}

\hypersetup{
   allcolors = {blue}
}

\begin{document}
\title{Efficient numerical method for reliable upper and lower bounds on homogenized parameters}
\author{Liya Gaynutdinova}\address{Department of Mathematics, Faculty of Civil Engineering, Czech Technical University in Prague, Thákurova 2077/7, Prague, Czech Republic}
\author{Martin Ladecký}\sameaddress{1}
\secondaddress{Department of Mechanics, Faculty of Civil Engineering, Czech Technical University in Prague, Thákurova 2077/7, Prague, Czech Republic}
\author{Aleš Nekvinda}\sameaddress{1}
\author{Ivana Pultarová}\sameaddress{1}
\author{Jan Zeman}\sameaddress{2}\secondaddress{Corresponding author: \href
{mailto:jan.zeman@cvut.cz}{jan.zeman@cvut.cz}}

\date{\today}

\begin{abstract}
A numerical procedure providing guaranteed two-sided bounds on the effective
coefficients of elliptic partial differential operators is presented. The upper bounds are obtained in a standard manner through the variational formulation of the problem and by applying the finite element method. To obtain the lower bounds we formulate the dual variational problem and introduce appropriate approximation spaces employing the finite element method as well. We deal with the 3D setting, which has been rarely considered in the literature so far. The theoretical justification of the procedure is presented and supported with illustrative examples.
\end{abstract}

\subjclass{35B27, 65N15, 49N15, 65F10}

\keywords{homogenization of linear elliptic PDEs; duality-based error estimates; finite element method; Helmholtz decomposition; conjugate gradient method}

\maketitle

\section{Introduction}

\baselineskip = 18pt

Composite materials are prominent in many areas of engineering because of their diverse useful properties. 
Such materials are often considered to posses a periodic structure. 
To obtain the effective properties of those materials, the homogenization theory~\cite{Jikov,Michel} was developed, and the recent decades {witnessed} a great progress of {efficient} numerical methods in homogenization.
They developed from the iterative fixed point 
computational schemes for Lippmann-Schwinger equation
using collocation principles~\cite{Moulinec1994,Moulinec1998} to the finite difference or variational approaches which use Krylov subspace solvers
and the fast Fourier {transform} 
algorithm to accelerate the convergence;
see e.g.~\cite{Brisard,Lucarini2019,Leuschner2018,Schneider2016,Schneider2017,ZemanVNM} 
or~\cite{Schneider2021,Lucarini,Gierden2022}
for recent comprehensive reviews. {As a result of these developments, the contemporary FFT-accelerated solvers are capable of up-scaling complex multi-physics processes defined on high-resolution voxel microstructures on conventional hardware. Yet, the issues of convergence and accuracy (or error estimates) of the resulting homogenized properties have gained attention only recently; see~\cite{Schneider2022nme,Ye:2023,Schneider:2023} for the current convergence results and~\cite{Ferrier:2023} for a recent study on error estimation.}

{In this work, we adopt} a duality{-based} framework, see e.g.~\cite{Briane,Serrano,
Wieckowski},
consider a 3D scalar elliptic partial differential equation with periodic data,
and introduce a numerical procedure
providing guaranteed and arbitrarily close 
upper and lower bounds on the homogenized data.
Our approach {rests on earlier results by Dvořák and Haslinger}~\cite{Dvorak,HaslingerDv}
where, however, only a 2D setting is considered; see also~\cite{Briane}. 
The main difficulty of {generalization to 3D 
consists in} choosing appropriate finite dimensional approximation spaces
to the periodic divergence-free solution space of the dual problems. 
The only paper dealing with a 3D setting is~\cite{VondrejcEtAl2015} 
where the approximation spaces are {formed by} smooth periodic 
spectral Fourier functions.  
There are two drawbacks of numerical computation with a truncated
Fourier expansion. Numerical solutions of 
differential equations with discontinuous data may 
suffer from ringing artifacts, e.g.~\cite[Section~3]{Leute2022} {for examples in the present context}, and computing the elements of
stiffness matrices exactly may become difficult {even for piecewise constant coefficients and grids fitted to coefficient jumps}~\cite{Monchiet2015,Vondrejc2016}.
However, computing elements of the linear system (almost) exactly is important for
providing reliable bounds on the homogenized data.

In this paper we consider a cuboid domain $Y\subset \mathbb{R}^3$
and a symmetric matrix function 
$A\in L^\infty(Y,{\mathbb R}^{3\times 3})$ which is 
uniformly positive definite a.e.~in $Y$, i.e.,
there exist $c_1,c_2>0$ such that $c_1(v,v)_{\mathbb{R}^3}\le 
(A(x)v,v)_{\mathbb{R}^3}\le c_2 (v,v)_{\mathbb{R}^3}$
for all $v\in\mathbb{R}^3$ and a.a.~$x\in Y$.
Here $(u,v)_{\mathbb{R}^3}$ denotes the inner product of
vectors $u$ and $v$ of the Euclidean space $\mathbb{R}^3$.
The related homogenization problem~\cite{Jikov} reads to find
a homogenized data matrix $A^*\in {\mathbb R}^{3\times 3}$ such that
\begin{equation}\label{min1}
(  A^*\xi,\xi)_{{\mathbb R}^3}=\inf_{u\in H^1_{\rm per}(Y,\mathbb{R})}\frac{1}{\vert Y\vert}
\int_Y(  A(\xi+\nabla u),\xi+\nabla u)_{{\mathbb R}^3}\, \dd x
\end{equation}
for all $\xi \in {\mathbb R}^3$, where $\vert Y\vert$ denotes the volume of $Y$. 
Our goal in this paper is to present a method that yields
guaranteed, feasible, and arbitrarily close 
upper and lower bounds on $A^*$.
Let $A^*_h$ be obtained as
a result of a minimization problem~\eqref{min1}
over only a finite-dimensional subspace 
$V_h\subset H^1_{\rm per}(Y,{\mathbb R})$; see Section~\ref{sec3} for details and definitions.
Matrix $A^*_h$ is an upper bound to $A^*$, 
$A^*\preccurlyeq A^*_h$, in the sense {that}
\begin{equation}\nonumber
(  A^*\xi,\xi)_{{\mathbb R}^3}\le (  A^*_h\xi,\xi)_{{\mathbb R}^3}\quad
\text{for all}\;\xi\in {\mathbb R}^3.
\end{equation}
A sufficiently close approximation $A^*_h$ to $A^*$ can be obtained by appropriate refining of the solution space $V_h$. 
It remains to obtain the computable sufficiently accurate guaranteed  
lower bounds for $A^*$. This sets the 
goal of this paper.

The outline of the paper is as follows.
In the next section we recall the Helmholtz decomposition 
of $L^2(Y,\mathbb{R}^3)$ into rotation-free and divergence-free 
$Y$-periodic functions. As the first result of this paper, we propose the approximation spaces for this decomposition
based on the finite element (FE) discretization using standard Lagrange finite elements. In Section~\ref{sec4} we suggest and theoretically prove a tool 
for obtaining the lower bounds on homogenized coefficients.
In contrast to~\cite{VondrejcEtAl2015,Dvorak}, we build the theoretical
justification on a simple optimization result which is proved here,
instead of on a rather involved perturbation-duality argument.
{At the end of the section, we propose a novel projection-based
approximate solution of the dual problem and its error estimate. This solution does not involve any iterative solution
process.}
The reader who searches for computational 
algorithms only may skip Sections~\ref{sec3} and~\ref{sec4}
and proceed {to} Section~\ref{sec5} directly { in which} we comment on the implementation of our procedure and present some numerical examples.
Some concluding remarks and the plans for further research are found in Section~\ref{sec6}.

\section{Helmholtz decomposition and approximation spaces}\label{sec3}

We consider a cuboid domain $Y\subset{\mathbb R}^3$, $Y=(0,a_1)\times (0,a_2) \times (0,a_3)$, $a_j>0$, $j=1,2,3$. Then $Y$ is Lipschitz, connected, and bounded.
We only deal with Euclidean space ${\mathbb R}^3$. 
Changing the setting to ${\mathbb R}^2$ is straightforward and fully described in~\cite{Dvorak}.

\subsection{Helmholtz decomposition of periodic fields}\label{sec222}

Let us recall for 
a domain $\Omega\subset\mathbb{R}^3$ the function space $L^2(\Omega,{\mathbb R}^k)$
of functions $u:\Omega\to\mathbb{R}^k$,
$u=(u_1,\dots,u_k)^T$,
with the inner product
$(u,v)=\int_\Omega (u,v)_{\mathbb{R}^k} \dd x$, and 
the function space $L^1(\Omega,{\mathbb R}^k)$
containing functions $u:\Omega\to {\mathbb R}^k$ 
with finite Lebesgue integrals 
$\int_\Omega \vert u_j(x)\vert\, \dd x$, $j=1,\dots,k$.

Considering the special case $\Omega=Y$,
the $Y$-periodic extension $u_{\rm per}$ of 
$u\in L^1(Y,{\mathbb R}^3)$ is defined as 
\begin{equation}\nonumber
u_{\rm per}(x)=u(x)\; \text{for a.a.}\; x\in Y,
\end{equation}
and
\begin{equation}\nonumber
u_{\rm per}(x)=u_{\rm per}(x+(a_1,0,0))
=u_{\rm per}(x+(0,a_2,0))=u_{\rm per}(x+(0,0,a_3))
\; \text{for a.a.}\;
x\in {\mathbb R}^3.
\end{equation}
The mean value $\langle u\rangle\in\mathbb{R}^k$
of $u\in L^1(Y,{\mathbb R}^k)$ is defined as
\begin{equation}\nonumber
\langle u\rangle=\frac{1}{\vert Y\vert}\int_Y u\, \dd x,
\end{equation}
where $\vert Y\vert=a_1a_2a_3$.
Any subspace of $L^1(Y,{\mathbb R}^3)$ of functions with zero mean 
is denoted by the subscript $_{{\rm mean},0}$.
Let us {introduce} the function spaces
\begin{eqnarray}\nonumber
H^{1}(Y,{\mathbb R}^k)&=&\{ u\in L^2(Y,{\mathbb R}^k);
\; \nabla u_j\in L^2(Y,{\mathbb R}^{3}),\; j=1,\dots,k\}\nonumber\\
H^{\rm curl}(Y,{\mathbb R}^3)&=&\{ u\in L^2(Y,{\mathbb R}^3);\; {\rm curl}\, u\in L^2(Y,{\mathbb R}^3)\}\nonumber\\
H^{\rm div}(Y,{\mathbb R}^3)&=&\{ u\in L^2(Y,{\mathbb R}^3);\; {\rm div}\, u\in L^2(Y,{\mathbb R})\}\nonumber
\end{eqnarray}
where $u:Y\to {\mathbb R}$ and $v:Y\to {\mathbb R}^3$ and where
the differential operators 
\begin{equation}\nonumber
\nabla u=\left(\frac{\partial}{\partial x_1},\frac{\partial}{\partial x_2},
\frac{\partial}{\partial x_3}\right)u,\quad
{\rm div}\, v=\nabla\cdot v,\quad {\rm curl}\, v=\nabla\times v
\end{equation}
are defined in the weak sense.
The functions with zero weak divergence in $Y$ are defined 
as
\begin{equation}\nonumber
H^{{\rm div},0}(Y,{\mathbb R}^3)=\{ u\in L^{2}(Y,{\mathbb R}^3);\; \int_Yu\cdot\nabla\phi\,{\rm d}x=0,\; \phi\in C^\infty_0(Y,\mathbb{R})\},
\end{equation}
see, e.g., \cite{Girault}.
Let the subscript $_{\rm loc}$ denote that a required property 
is fulfilled on all compact subsets of ${\mathbb R}^3$.
Let us define
\begin{eqnarray}\nonumber
H^1_{\rm per}(Y,{\mathbb R}{^k})
&=&
\{ u\in H^1(Y,{\mathbb R}{^k});\; 
u_{\rm per}\in H^1_{\rm loc}({\mathbb R}^3,{\mathbb R}{^k})\}
\nonumber
\\
H^{{\rm div},0}_{\rm per}(Y,{\mathbb R}^3)&=&\{ u\in H^{{\rm div},0}(Y,{\mathbb R}^3);\; u_{\rm per}\in H^{{\rm div},0}_{\rm loc}({\mathbb R}^3,{\mathbb R}^3)\}\nonumber\\
\nabla H^1_{\rm per}(Y,{\mathbb R})&=&
\{ u\in L^2(Y,{\mathbb R}^3);\; u=\nabla v,\, v\in  H^1_{\rm per}(Y,{\mathbb R})\}.\nonumber
\end{eqnarray}

From~\cite{Ranocha,Deriaz} or~\cite{Schweizer2017} and~\cite{Jikov}, 
we have the following theorem.
\begin{thrm}\label{THM1}
A vector field $u\in L^2(Y,{\mathbb R}^3)$ satisfies ${\rm div}\, u=0$
if and only if there exists a vector potential 
$v\in H^{\rm curl}(Y,{\mathbb R}^3)$ for which $u={\rm curl}\, v$.
\end{thrm}
In our setting, we will decompose
$L^2(Y,\mathbb{R}^3)$ into periodic irrotational (potential)
fields ($\nabla H^1_{\rm per}
(Y,\mathbb{R})$) and periodic divergence free ($ H^{{\rm div},0}_{\rm per}(Y,\mathbb{R}^3)$) fields.
Let us first recall the following divergence theorem~\cite{Brezzi,Girault}: 
\begin{lmm}\label{Pok1}
There exists a continuous linear operator $\gamma:H^{{\rm div},0}(Y,\mathbb{R}^3)\to H^{-1/2}(\partial Y,\mathbb{R})$ such that
\begin{equation}\nonumber
\gamma u=u\cdot n|_{\partial Y}\quad \text{for}\quad u\in C^\infty(\overline{Y},\mathbb{R}^3).
\end{equation}
For $u\in H^{{\rm div},0}(Y,\mathbb{R}^3)$ and $v\in H^1(Y,\mathbb{R})$
we have
\begin{equation}\label{divthm}
\int_Yu\cdot \nabla v\, \dd x=\langle \gamma u, Tv\rangle, 
\end{equation}
where $Tv$ is the trace of $v$ on the boundary $\partial Y$ and $\langle\phi,
\psi\rangle$ is a duality pairing between
$\psi\in H^{1/2}(\partial Y,\mathbb{R})$ and 
$\phi\in H^{-1/2}(\partial Y,\mathbb{R})$.
\end{lmm} 

\begin{rmrk}
Images of $\gamma$ defined by Lemma~\ref{Pok1} can be denoted by
$\gamma u=u\cdot n$ and called
normal fluxes of $u\in H^{{\rm div},0}(Y,\mathbb{R}^3)$ on $\partial Y$.
\end{rmrk}

{ It is well known, see e.g.,~\cite{Girault}, that} 
for every function $u\in L^2(Y,\mathbb{R}^3)$ there exist
functions $v$ and $w$ such that $u=\nabla v+{\rm curl}\,w$ where
$v\in H^1(Y,\mathbb{R})$,  $w\in H^1(Y,\mathbb{R}^3)$,
and $\gamma(u-\nabla v)=0$.
In our approach, however, we will deal with $Y$-periodic functions only,
which are less frequent in the literature. 
{Although the main results needed later are stated in~\cite{Jikov}, we present their detailed derivation 
in Appendix~\ref{apxA} for the readers convenience.}

\begin{thrm}\label{charU}
Let us define 
\begin{equation}\nonumber
W=\{u\in L^2(Y,{\mathbb R}^3);\, \int_Yu\cdot \nabla \phi\, dx=0
\;\; {\text{for all}}\;\; \phi\in H^1_{\rm per}(Y,{\mathbb R})\}.
\end{equation}
Then
\begin{equation}\nonumber
W=H^{{\rm div},0}_{\rm per}(Y,{\mathbb R}^3).
\end{equation}
\end{thrm}

\begin{rmrk}\label{rem111}
The spaces 
$H^{{\rm div},0}_{\rm per}(Y,{\mathbb R}^3)$,
$\nabla H^1_{\rm per}(Y,{\mathbb R})$, and 
$H^{{\rm div},0}_{\rm per,mean,0}(Y,{\mathbb R}^3)$
are closed subspaces of $L^2(Y,{\mathbb R}^3)$;
the proof proceeds in the standard way using the Green theorem and 
the Cauchy-Schwartz inequality, see e.g.~\cite{Girault}.
Then from Theorem~\ref{charU} it follows that there is an $L^2(Y;\mathbb{R}^3)$-orthogonal decomposition ({the} Helmholtz decomposition)
\begin{equation}\label{decomp}
L^2(Y,{\mathbb R}^3)=H^{{\rm div},0}_{\rm per}(Y,{\mathbb R}^3)
\oplus \nabla H^1_{\rm per}(Y,{\mathbb R}).
\end{equation}
Moreover, we have
\begin{equation}\nonumber
L^2(Y,{\mathbb R}^3)=H^{{\rm div},0}_{\rm per,mean,0}(Y,{\mathbb R}^3)
\oplus \nabla H^1_{\rm per}(Y,{\mathbb R})\oplus \mathbb{R}^3.
\end{equation}
\end{rmrk}

\subsection{Approximation spaces}

As we will see in Section~\ref{sec4}, for the numerical computation of the upper and lower bounds on $A^*$ it is necessary to accurately approximate 
functions from $\nabla H^1_{\rm per}(Y,{\mathbb R})$ and from $H^{{\rm div},0}_{\rm per}(Y,{\mathbb R}^3)$, respectively.
A common approach to approximate functions of $\nabla H^1_{\rm per}(Y,{\mathbb R})$ is to use the FE basis functions, also called shape functions and their derivatives.
Approximation of functions of $H^{{\rm div},0}_{\rm per}(Y,{\mathbb R}^3)$ is 
less common. In the literature, they are also called  
statically admissible functions, e.g.~\cite{Krizek1983}, and special unisolvent FE spaces were proposed, e.g.,~\cite{HH1976,Nedelec1980}. 
In the remaining part of this section, we show how standard FE basis functions can be used to approximate 
$H^{{\rm div},0}_{\rm per,mean,0}(Y,{\mathbb R}^3)$. 
The absence of unisolvence is compensated by an easy implementation. 

{To this goal, let us define the following} matrices
\begin{equation}\label{Q123}
Q_1=\left(\begin{array}{ccc}
0&-1&0\\
1&0&0\\
0&0&0
\end{array}
\right),\quad
Q_2=\left(\begin{array}{ccc}
0&0& 1\\
0&0&0\\
-1&0&0
\end{array}
\right)
\quad{\text{and}}
\quad
Q_3=\left(\begin{array}{ccc}
0&0&0\\
0&0&-1\\
0&1&0
\end{array}
\right).
\end{equation}

\begin{lmm}\label{Lem11}
The closure of 
\begin{equation}\label{WQ}
W_Q={\rm span}\,\left(Q_1\nabla H^1_{\rm per}(Y,{\mathbb R})
\cup Q_2\nabla H^1_{\rm per}(Y,{\mathbb R})
{\cup Q_3\nabla H^1_{\rm per}(Y,{\mathbb R})}\right)
\end{equation}
with respect to the metric of $L^2(Y,{\mathbb R}^3)$  
equals to $H^{{\rm div},0}_{{\rm per},{\rm mean},0}(Y,{\mathbb R}^3)$.\end{lmm}
\begin{proof}
Let $u\in W_Q$. Then there exist potentials
$\psi_1$, $\psi_2$ { and $\psi_3$} in $H^1_{\rm per}(Y,{\mathbb R})$ such that $u=Q_1\nabla\psi_1+Q_2\nabla\psi_2{+Q_3\nabla\psi_3}$.
Then for $\phi\in C^\infty_{\rm per}(Y,{\mathbb R})$
\begin{equation}\nonumber
\int_Y (Q_j\nabla\psi_j)\cdot\nabla\phi\, \dd x
=-\int_Y \psi_j\cdot {\rm div}\,Q_j^T\nabla\phi\, \dd x+\int_{\partial Y}(n\cdot Q_j^T\nabla \phi)\psi_j\, \dd s=0,\quad j=1,2{,3,}
\end{equation} 
because ${\rm div}\,Q_j^T\nabla\phi=0$, and $(n\cdot Q_j^T\nabla \phi)\psi_j$ has opposite values on the
opposite faces of $Y$. Then
\begin{equation}\nonumber
\int_Yu\cdot\nabla\phi\, \dd x=\int_Y (Q_1\nabla\psi_1+Q_2\nabla\psi_2{+Q_3\nabla\psi_3})\cdot\nabla\phi\, \dd x=0.
\end{equation} 
Since $C^\infty_{\rm per}(Y,{\mathbb R})$ is dense in $H^1_{\rm per}(Y,\mathbb{R})$, we get $u\in H^{\rm div,0}_{\rm per}(Y,\mathbb{R}^3)$ by~Theorem~\ref{charU}.
Since $\psi_j$ is periodic, $Q_j\nabla\psi_j$ is periodic as well and
the mean of each component is zero.
This yields that $W_Q$
is a subspace of $H^{{\rm div},0}_{{\rm per},{\rm mean},0}(Y,{\mathbb R}^3)$.
{ Since $H^{{\rm div},0}_{{\rm per},{\rm mean},0}(Y,{\mathbb R}^3)$ is closed, then it also contains the closure of $W_Q$.}

{Now consider $u\in H^{{\rm div},0}_{{\rm per},{\rm mean},0}(Y,{\mathbb R}^3)$. Because $H^{{\rm div},0}_{{\rm per},{\rm mean},0}(Y,{\mathbb R}^3) \subset L^2(Y,{\mathbb R}^3)$, the field $u$ 
admits the Fourier series representation}
\begin{equation}\label{uF}
u(x)=\sum_{j\in{\mathbb Z}^3{,\, j\ne 0}}u_{j}(x), \quad u_{j}(x)=c_j
{\rm exp}\,(\ii j\cdot x),\quad c_j \in {\mathbb C}^3,
\end{equation}
{where we assumed, for notation simplicity, $Y=(0,2\pi)^3$.} Note that 
\begin{equation}\nonumber
{\rm div}\,u_{j}=\ii\, c_j\cdot j\; {\rm exp}\,(\ii j\cdot x)
\end{equation}
where $j=(j_1,j_2,j_3)^T\in{\mathbb Z}^3$, $x=(x_1,x_2,x_3)^T\in Y$.
Every term $u_{j}$ of the expansion~\eqref{uF} is $Y$-periodic. 
Using ${\rm div}\,u=0$ and Theorem~\ref{charU} with $\phi={\rm exp}\,
(\ii j\cdot x)$, we get {$j\cdot c_j=0$ for all $j \in {\mathbb Z}^3$, which implies that that each Fourier coefficient $c_j$ must be of a from $c_j = t_j \times j$ with $t_j \in \mathbb{C}^3$. Therefore,}
\begin{equation}\nonumber
c_j 
= 
{
(t_{j,1}Q_1+t_{j,2}Q_2 + t_{j,3}Q_3 ) j 
},\quad t_{j,1}, t_{j,2}{, t_{j,3}  }\in{\mathbb C},
\end{equation}
{and} we can write
\begin{equation}\label{FT1}
u_{j}(x)=(t_{j,1}Q_1+t_{j,2}Q_2{+t_{j,3}Q_3})\nabla {\rm exp}\,(\ii j\cdot x),\quad t_{j,1},t_{j,2}{,t_{j,3}}\in{\mathbb C}.
\end{equation} 
Let us consider a sufficiently close (in the norm of $L^2(Y,{\mathbb R}^3)$) approximation $u^{(m)}$ to $u$ with a finite number of terms
\begin{equation}\label{uFfinite}
u^{(m)}(x)=\sum_{\vert j\vert\le m}u_j(x)
=\sum_{\vert j\vert\le m}c_j{\rm exp}\,(\ii j\cdot x).
\end{equation}
Then
\begin{equation}\nonumber
u^{(m)}(x)=\sum_{\vert j\vert\le m}c_j{\rm exp}\,(\ii j\cdot x)=
\sum_{\vert j\vert\le m}\left(t_{j,1}Q_1+t_{j,2}Q_2
{+t_{j,3}Q_3}\right)
\nabla {\rm exp}\,(\ii j \cdot x)
=Q_1\nabla\psi_1(x)+Q_2\nabla\psi_2(x){+Q_3\nabla\psi_3(x)}
\end{equation}
where 
\begin{equation}\nonumber
\psi_s(x)=\sum_{\vert j\vert\le m}t_{j,s}\,{\rm exp}\,(\ii j\cdot x),
\end{equation}
$\psi_s\in C^\infty_{\rm per}(Y,{\mathbb R})$,
$s=1,2,3$. 
Such a $u^{(m)}\in W_Q$ can be arbitrarily close to $u$.

{Considering general cell} $Y=(0,a_1)\times (0,a_2)\times (0,a_3)$ yields  
\begin{equation}\nonumber
{\rm div}\,u_{j}=\frac{\ii}{2\pi}\sum_{k=1}^3\,(c_j)_k a_kj_k\; {\rm exp}\,\left(\frac{\ii}{2\pi} \sum_{k=1}^3 {a_k} j_k x_k\right).
\end{equation}
The coefficients of the field $u$ with zero divergence are then
\begin{equation}\nonumber
c_j 
{
=
(t_{j,1}Q_1+t_{j,2}Q_2 + t_{j,3}Q_3 )\, a \odot j  
}
\end{equation}
{where $a \odot j$ denotes the element-wise (Hadamard) product.} 
Thus~\eqref{FT1} remains valid, and the statement 
of the lemma as well.
\end{proof}

\begin{crllr}\label{coro1}
From Lemma~\ref{Lem11} it follows that the following orthogonal decomposition holds
\begin{equation}\label{decompos1}
L^2(Y,{\mathbb R}^3)=\nabla H^1_{\rm per}(Y,{\mathbb R})
\oplus \overline{W_Q}\oplus {\mathbb R}^3,
\end{equation}
where 
\begin{equation}\label{WQ3}
W_Q={\rm span}\,\left(Q_1\nabla H^1_{\rm per}(Y,{\mathbb R})
\cup Q_2\nabla H^1_{\rm per}(Y,{\mathbb R})
\cup Q_3\nabla H^1_{\rm per}(Y,{\mathbb R})\right).
\end{equation}
{In analogy to to the notation $\nabla H^1_{\rm per}(Y,\mathbb{R})$,
we can use the notation $W_Q={\rm curl}\, H^1_{\rm per}(Y,\mathbb{R}^3)$ and write 
\begin{equation}
    L^2(Y,{\mathbb R}^3)=\nabla H^1_{\rm per}(Y,{\mathbb R})
\oplus {\rm curl}\, H^1_{\rm per}(Y,{\mathbb R}^3)   \oplus {\mathbb R}^3.
\end{equation}
}
\end{crllr}

\begin{rmrk}
From density of standard FE spaces in $H^1_{\rm per}(Y,\mathbb{R})$ (for example, continuous piece-wise linear periodic functions; see e.g.~\cite{Brenner})
and from~\eqref{decompos1}, we obtain that the set of
gradients of the FE basis 
functions is dense in $\nabla H^1_{\rm per}(Y,\mathbb{R})$
in the $L^2(Y,\mathbb{R}^3)$ metric{, and the set of the same gradients multiplied by matrices $Q_j$, $j=1,2,3$, is dense in ${\rm curl}\, H^1_{\rm per}(Y,{\mathbb R}^3)$ in the same metric.} An interested reader may compare our approach with a related study~\cite{Krizek1986} for non-periodic domains and different boundary conditions. {
As such, the presented approach extends and combines the two previous frameworks proposed in~\cite{Dvorak,VondrejcEtAl2015}. Specifically, our characterization of the space ${\rm curl}\, H^1_{\rm per}(Y,{\mathbb R}^3)$ generalizes the two-dimensional stream function technique of Dvo\v{r}\'{a}k~\cite{Dvorak} to three dimensions. In~\cite{VondrejcEtAl2015}, the same space was constructed from $L^2(Y,{\mathbb R}^3)$ using a projection operator with an explicit expression in the Fourier space. Note that the dimensions of the spaces spanned with the trigonometric and finite element bases coincide in both two and three dimensions.}
\end{rmrk}

\section{Upper and lower bounds on homogenized coefficients}\label{sec4}

In this section, we focus on obtaining the upper and lower bounds on $A^*$,
where $A^*$ is defined by~\eqref{min1}.
By bounds, we mean symmetric positive definite matrices
{$C^*_1, C^*_2\in\mathbb{R}^{3\times 3}$} such that 
\begin{equation}\label{upper}
({C^*_1}\xi,\xi)_{\mathbb{R}^3}\le (A^*\xi,\xi)_{\mathbb{R}^3}\le 
({C^*_2}\xi,\xi)_{\mathbb{R}^3},\quad \text{for all}\;\xi\in\mathbb{R}^3,
\end{equation}
which we denote by ${C^*_1}\preccurlyeq A^*\preccurlyeq {C^*_2}$.

Recall that
$A:Y\to\mathbb{R}^{3\times 3}$ is 
uniformly positive definite and bounded a.e.~in $Y$.
{It is well known that the minimum of quadratic functional defined
by~\eqref{min1}} is attained for $u=u_\xi$
satisfying 
\begin{equation}\label{form111}
    \int_Y(A\nabla u_\xi,\nabla v)_{\mathbb{R}^3}
    \, \dd x=-
    \int_Y(A\xi,\nabla v)_{\mathbb{R}^3}
    \, \dd x,\qquad \text{for all}\; v\in H^1_{\rm per}(Y;\mathbb{R}).
\end{equation}
Observe that the mapping $\xi\to u_\xi$ is linear, which yields for all $\xi\in\mathbb{R}^3$,
\begin{eqnarray}
(A^*\xi,\mu)_{\mathbb{R}^3}&=&
\frac{1}{2}\left(
(A^*(\xi+\mu),\xi+\mu)_{\mathbb{R}^3}-
(A^*\xi,\xi)_{\mathbb{R}^3}-
(A^*\mu,\mu)_{\mathbb{R}^3}
\right)\nonumber\\
&=&\frac{1}{\vert Y\vert}\int_Y
(A(\xi+\nabla u_\xi),\mu+\nabla u_\mu)_{\mathbb{R}^3}\,\dd x.\nonumber
\end{eqnarray}
By~\eqref{form111} we also have
\begin{equation}\label{min11}
(A^*\xi,\mu)_{\mathbb{R}^3}=
\frac{1}{\vert Y\vert}\int_Y
(A(\xi+\nabla u_\xi),\mu)_{\mathbb{R}^3}\,\dd x=
\frac{1}{\vert Y\vert}\int_Y
(A\xi,\mu+\nabla u_\mu)_{\mathbb{R}^3}\,\dd x.
\end{equation}

\subsection{Upper bounds}

The minimizing solution $u_\xi\in H^1_{\rm per}(Y,\mathbb{R})$ in~\eqref{min1} is obtained by solving~\eqref{form111}.
A numerical solution can be found on a finite-dimensional subspace 
$U_h$ of $H^1_{\rm per}(Y,\mathbb{R})$, usually generated by FE functions. 
Since~\eqref{min1} is formulated as a variational problem where the minimum is to be reached over some infinite dimensional space, any 
minimum {over} a subspace $U_h\subset H^1_{\rm per}(Y,{\mathbb R})$ or even any functions from
$H^1_{\rm per}(Y,\mathbb{R})$ yield an upper bound to $A^*$ in the sense~\eqref{upper}. 
More precisely, let $\widetilde{u}^{(j)}\in U_h$ be arbitrary and
$\xi=e^{(j)}$, $j=1,2,3$,
where $e^{(j)}$ is the $j$-th column of $3\times 3$ identity matrix.
Let 
\begin{equation}\nonumber
(A_h^*)_{jk}=\frac{1}{\vert Y\vert}
\int_Y(  A(e^{(k)}+\nabla \widetilde{u}^{(k)}),e^{(j)}+\nabla \widetilde{u}^{(j)})_{{\mathbb R}^3}\, \dd x,
\quad j,k=1,2,3.
\end{equation}
Then $A^*_h$ is an upper bound to $A^*$ in the sense of~\eqref{upper}, $A^*\preccurlyeq A^*_h$.
A sufficiently close approximation can be obtained by refining the solution space $U_h$ { and thus by
having $u^{(j)}$} closer to their minimizing counterparts in $H^1_{\rm per}(Y,\mathbb{R})$ in general. 
It remains to obtain the computable 
sufficiently close guaranteed lower bounds, which is the
goal of the { next two subsections}.

\subsection{Lower bounds}

Our estimation method relies on the perturbation duality theorem.
Various forms of the theorem can be found in the literature;
see e.g.~{\cite{EkelandTemam,HH1976,Jikov}}. 
To keep the paper self-contained, we present and prove a certain form of the duality theorem as Lemma~\ref{lem2}, which is significantly simpler 
than the version in~\cite{EkelandTemam} used in, e.g.,
\cite{Dvorak,VondrejcEtAl2015}.
Using Lemma~\ref{lem2}, the problem~\eqref{min1} is then 
transformed to the dual one; see Lemma~\ref{Lem1}.
{Both lemmas and their proofs are found in Appendix~\ref{apxB}.}

We now formulate the duality principle for the homogenization problems.
In the rest of this {paper}, we will use the following notation
\begin{align*}\nonumber
H & = L^2(Y,{\mathbb R})^3, 
&
H_0 & = \{u\in L^2(Y,{\mathbb R}^3),\, \langle u\rangle = 0\},
&
U &= H^1_{\rm per}(Y,{\mathbb R}),
\\
V & = \nabla U,
&
W & = H^{{\rm div},0}_{\rm per}(Y,{\mathbb R}^3),
&
W_0 &= H^{{\rm div},0}_{{\rm per},{\rm mean}, 0}(Y,\mathbb{R}^3).\nonumber
\end{align*}
Note that the mean $\langle v\rangle=0$ for any $v\in V$.
According to Remark~\ref{rem111} we have 
\begin{equation}\nonumber
H=V\oplus W\qquad \text{and}\qquad H_0=V\oplus W_0.
\end{equation}
Since $V$ is a closed subspace in $H$, we can use 
minimum instead of infimum in~\eqref{min1}.

Let us define {\it the dual problem} to~\eqref{min1}:
Find a matrix $B^*\in{\mathbb R}^{3\times 3}$ such that for every $\alpha\in {\mathbb R}^3$
\begin{equation}\label{min2}
(B^*\alpha,\alpha)_{{\mathbb R}^3}=\frac{1}{\vert Y\vert}
\min_{w\in W_0}
\int_Y(  A^{-1}(\alpha+w),
\alpha+w)_{{\mathbb R}^3}\, \dd x.
\end{equation} 
Note that $B^*$ is symmetric and positive definite.

\begin{thrm}\label{lem55}
Let $B^*$ be defined by~\eqref{min2}. Then we have
\begin{equation}\nonumber
(B^*)^{-1}=A^*.
\end{equation}
Moreover, let $\alpha=A^*\xi$,  
$u_\xi\in U$ and $w_\alpha\in W_0$ be the minimizers in~\eqref{min1} and~\eqref{min2}, respectively, and let $w_1$ be given by
\begin{equation}\label{AEj}
A^*\xi+w_1=A(\xi+\nabla u_\xi),\quad \xi\in {\mathbb R}^3.
\end{equation}
Then $w_1=w_\alpha$.
\end{thrm}
Following {Theorem~\ref{lem55}} we can get $A^*=(B^*)^{-1}$ by solving
the dual problem~\eqref{min2}.
Since the dual problem is a variational minimization problem,
for any approximate solution $B^*_h$ we obtain an upper bound to 
$B^*$ which provides a lower bound to $A^*$,
$(B^*_h)^{-1}\preccurlyeq (B^*)^{-1}\preccurlyeq A^*$ {, e.g., \cite[Corollary 7.7.4]{horn_johnson_2012}}. 
For completeness, we {present this well-known equivalence:}
\begin{lmm}
Let $A,B\in{\mathbb R}^{n\times n}$ be symmetric positive definite
matrices. Then 
\begin{equation}\label{ekviv1}
v^TAv\le v^TBv,\quad \text{ for all } v\in {\mathbb R}^n\quad \iff
\quad
v^TB^{-1}v\le v^TA^{-1}v,\quad \text{ for all } v\in {\mathbb R}^n.
\end{equation}
\end{lmm}

{

\subsection{Numerical solution}\label{secnumsol}

Numerical solution of the minimizers $u_\xi\in U$ and 
$w_\alpha\in W_0$ of~\eqref{min1} and~\eqref{min2},
respectively, usually yields some inexact $\widetilde{u}_\xi\in U$ and 
$\widetilde{w}_\alpha\in W_0$, respectively. Then obviously
\begin{eqnarray}
(A^*\xi,\xi)_{\mathbb{R}^3}&\le&
\frac{1}{\vert Y\vert}
\int_Y(  A(\xi+\widetilde{u}_\xi),
\xi+\nabla\widetilde{u}_\xi)_{{\mathbb R}^3}\, \dd x\nonumber\\
(B^*\alpha,\alpha)_{\mathbb{R}^3}&\le&
\frac{1}{\vert Y\vert}
\int_Y(  A^{-1}(\alpha+\widetilde{w}_\alpha),
\alpha+\widetilde{w}_\alpha)_{{\mathbb R}^3}\, \dd x.\nonumber
\end{eqnarray}
The difference between the left and right-hand sides is handled
in the following lemma; see also~\cite{Schneider2022nme,Ye:2023}, in which
$\lambda_1(A)$ and $\lambda_2(A)$ denote 
the essential infimum and the essential supremum of minimal and maximal eigenvalues of $A$ over $Y$, respectively.

\begin{lmm}\label{lem_dist}
Let $\widetilde{u}_\xi,\widetilde{u}_\nu\in U$ and $\widetilde{w}_\alpha,
\widetilde{w}_\beta\in W_0$ be some
approximations to $u_\xi,u_\nu$ and $w_\alpha,w_\beta$, respectively,
for some $\xi,\nu,\alpha,\beta\in\mathbb{R}^3$. Then
\begin{eqnarray}
&&\left\vert
\int_Y(  A(\xi+\nabla\widetilde{u}_\xi),
\nu+\nabla\widetilde{u}_\nu)_{{\mathbb R}^3}\, \dd x-\vert Y\vert
(A^*\xi,\nu)_{\mathbb{R}^3}\right\vert\nonumber\\
&&=
\left\vert \int_Y(  A(\nabla\widetilde{u}_\xi-\nabla u_\xi),
\nabla\widetilde{u}_\nu-\nabla u_\nu)_{{\mathbb R}^3}\, 
\dd x\right\vert\le
\lambda_2(A)
\Vert \widetilde{u}_\xi-u_\xi\Vert_U  
\Vert \widetilde{u}_\nu-u_\nu\Vert_U\nonumber\\
&&\left\vert
\int_Y(  A^{-1}(\alpha+\widetilde{w}_\alpha),
\beta+\widetilde{w}_\beta)_{{\mathbb R}^3}\, \dd x-\vert Y\vert
(B^*\alpha,\beta)_{\mathbb{R}^3}\right\vert\nonumber\\
&&=
\left\vert
\int_Y(  A^{-1}(\widetilde{w}_\alpha-w_\alpha),
\widetilde{w}_\beta-w_\beta)_{{\mathbb R}^3}\, \dd x\right\vert\nonumber\le
\frac{1}{\lambda_1(A)}
\Vert \widetilde{w}_\alpha-w_\alpha\Vert_{H}
\Vert \widetilde{w}_\beta-w_\beta\Vert_{H}.
\nonumber
\end{eqnarray}
\end{lmm}
\begin{proof}
We first deal with the primal problem. Note that for all 
$\xi\in \mathbb{R}^3$ and $v\in U$ we have
\begin{equation}\nonumber
  \int_Y\left( A(\xi+\nabla u_\xi),\nabla v\right)\, \dd x=0. 
\end{equation}
Then
\begin{eqnarray}
   && \int_Y(  A(\xi+\nabla\widetilde{u}_\xi),
\nu+\nabla\widetilde{u}_\nu)_{{\mathbb R}^3}\, \dd x-
\vert Y\vert(A^*\xi,\nu)_{\mathbb{R}^3}\nonumber\\
&&=  \int_Y(  A(\xi+\nabla\widetilde{u}_\xi),
\nu+\nabla\widetilde{u}_\nu)_{{\mathbb R}^3}\,\dd x-\int_Y
 (A(\xi+\nabla{u}_\xi),
\nu+\nabla{u}_\nu)_{{\mathbb R}^3}
\, \dd x \nonumber\\
&&=\int_Y  
(A\xi,\nabla\widetilde{u}_\nu)_{{\mathbb R}^3}+
(A\nabla\widetilde{u}_\xi,\nu)_{{\mathbb R}^3}+
(A\nabla\widetilde{u}_\xi,\nabla\widetilde{u}_\nu)_{{\mathbb R}^3}-
(A\nabla{u}_\xi,\nu)_{{\mathbb R}^3}
\, \dd x \nonumber\\
&&=\int_Y  
-(A\nabla u_\xi,\nabla\widetilde{u}_\nu)_{{\mathbb R}^3}
-(A\nabla\widetilde{u}_\xi,\nabla u_\nu)_{{\mathbb R}^3}+
(A\nabla\widetilde{u}_\xi,\nabla\widetilde{u}_\nu)_{{\mathbb R}^3}+
(A\nabla{u}_\xi,\nabla u_\nu)_{{\mathbb R}^3}
\, \dd x \nonumber\\
&&=
\int_Y\left( A(\nabla\widetilde{u}_\xi-\nabla u_\xi),\nabla\widetilde{u}_\nu-\nabla u_\nu\right)_{{\mathbb R}^3}\, \dd x.\nonumber
\end{eqnarray}
Absolute value of the last term is obviously bounded from above by $\lambda_2(A)\Vert \widetilde{u}_\xi-u_\xi\Vert_U
\Vert \widetilde{u}_\nu-u_\nu\Vert_U$. The proof for the dual form is analogous.
\end{proof}

Let us now assume that we use some finite-dimensional discretization space $U_h\subset U$.
For $\xi\in \mathbb{R}^3$, $\Vert\xi\Vert_{\mathbb{R}^3}=1$,
we obtain $\widetilde{u}_\xi\approx u_\xi$,
$\widetilde{u}_\xi\in U_h$. 
Let us assume that there exists $\delta_1>0$, such that 
$\Vert \widetilde{u}_\xi-u_\xi\Vert_U\le \delta_1$.
Similarly, for two other vectors which together with $\xi$
form an orthonormal basis of $\mathbb{R}^3$,
let us obtain approximations of all elements of $A^*$;
let us call such a matrix $\widetilde{A}^*$ and assume that
\begin{equation}\nonumber
\Vert \widetilde{A}^*-A^*\Vert_{2} \le \delta_2.
\end{equation}
We now consider some discretization space $W_{0,h}\subset W_0$
and want to use~\eqref{AEj} to obtain some $\widetilde{w}\in W_{0,h}$ which would be relatively close to $w_{\widetilde\alpha}$, where
$\widetilde{\alpha}=\widetilde{A}^*\xi$.
Let us denote 
\begin{equation}\nonumber
w_H=A(\xi+\nabla \widetilde{u}_\xi)-\widetilde{A}^*\xi.
\end{equation}
The field $w_H\in H$ may not belong to $W_{0,h}$ (and even to $W_0$),
because FE spaces, in general, 
are not endowed with the discrete Helmholtz decomposition,
unlike, for example, discrete Fourier bases; see, e.g.~\cite{VondrejcEtAl2015}.
Therefore, let $\widetilde{w}\in W_{0,h}$ be obtained as the $H$-orthogonal projection 
of $w_H$ onto $W_{0,h}$
\begin{equation}\nonumber
(\widetilde{w},v)_H=(w_H,v)_H,\quad v\in W_{0,h}.
\end{equation}
If $Y$ is a cuboid and the grid is regular, this projection would be easy to obtain numerically using the fast discrete Fourier
transform; see, e.g.~\cite{LadeckyEtAl22}.
Then the question remains what is the difference between
$\widetilde{w}$ and $w_{\widetilde\alpha}\in W_0$, the minimizer of~\eqref{min2} with $\alpha:=\widetilde\alpha$.
According to Lemma~\ref{lem_dist}, we would then immediately obtain the
bound to the error of computing  $(B^*\widetilde{\alpha},\widetilde{\alpha})_{\mathbb{R}^3}$
by using $\widetilde{w}$ instead of the exact minimizer $w_{\widetilde\alpha}$.
As we will see in the next lemma, the rate of convergence for the lower bound of $A^*$ obtained in this way is the same as for the upper bound.
\begin{lmm}\label{lem:projection_error}
Let $\xi\in \mathbb{R}^3$, $\Vert \xi\Vert_2=1$, be arbitrary.
Let $u_\xi$ be the minimizer in~\eqref{min1} and let 
$\widetilde {u}_\xi\in U$ be such that 
$\Vert \widetilde{u}_\xi - u_\xi\Vert_U<\delta_1$.
Let $\widetilde{A}^*
\in\mathbb{R}^{3\times 3}$ be an approximation to $A^*$
such that $\Vert\widetilde{A}^* - A^*\Vert_2<\delta_2$. 
Denote $\widetilde{\alpha}=\widetilde{A}^*\xi$ and 
\begin{equation}\nonumber
w_H=A\xi+A\nabla\widetilde{u}_\xi-\widetilde{A}^*\xi.
\end{equation}
Let $\widetilde{w}\in W_{0,h}$ be an $H$-orthogonal
projection of $w_H$ onto $W_{0,h}$, i.e.
\begin{equation}\nonumber
(\widetilde{w},v)_H=(w_H,v)_H,\quad v\in W_{0,h}.
\end{equation}
Then 
\begin{equation}
\Vert \widetilde{w}-w_{\widetilde{\alpha}}\Vert_H\le c_1\delta_1+c_2\delta_2+c  h,
\end{equation}
where $w_{\widetilde{\alpha}}\in W_0$ is the exact minimizer in~\eqref{min2}
for $\alpha:=\widetilde{\alpha}$, and $c_1,c_2,c>0$ depend on 
$\vert Y\vert$ and on essentially extremal eigenvalues of $A^*$ and $A$ on $Y$.
\end{lmm}
\begin{proof}
We first estimate  $\Vert w_H-\widetilde{w}\Vert_H$.
Let us denote $d_H=w_H-\widetilde{w}$.
Since $d_H\in H$, due to 
the Helmholtz decomposition~\eqref{decomp}
(see also Corollary~\ref{coro1})
we can decompose it into three pair-wise
$H$-orthogonal parts,
$d_H=d_W+d_V+d_C$ where $d_W\in W_0$, $d_V\in V$ and $d_C\in \mathbb{R}^3$. 
Therefore we have 
$\Vert d_H\Vert_H^2=\Vert d_W\Vert_H^2+\Vert d_V\Vert_H^2+
\Vert d_C\Vert_H^2$. 
Using~\eqref{form111} we have for all $v\in V$
\begin{eqnarray}\nonumber
 (d_V,v)_H&=&(w_H-\widetilde{w},v)_H=
(A\xi+A\nabla\widetilde{u}_\xi-\widetilde{A}^*\xi-\widetilde{w},v)_H=
(A\xi+A\nabla\widetilde{u}_\xi,v)_H\nonumber\\
&=&(A\xi+A\nabla\widetilde{u}_\xi-A\nabla{u}_\xi+A\nabla{u}_\xi,v)_H
=(A(\nabla\widetilde{u}_\xi-\nabla{u}_\xi),v)_H.\nonumber
\end{eqnarray}
Then
\begin{equation}\nonumber
 \vert (d_V,v)_H\vert \le \lambda_2(A)\Vert 
\nabla\widetilde{u}_\xi-\nabla u_\xi\Vert_H \Vert v\Vert_H\le \lambda_2(A)\delta_1
\Vert v\Vert_H,
\end{equation}
and thus $\Vert d_V\Vert_H\le \lambda_2(A)\delta_1$.
We now estimate $\Vert d_C\Vert_H$. 
Using~\eqref{min11}, we have for all $\nu\in\mathbb{R}^3$
\begin{eqnarray}\nonumber
 (d_C,\nu )_H&=&(w_H-\widetilde{w},\nu)_H=
(A\xi+A\nabla\widetilde{u}_\xi-\widetilde{A}^*\xi-\widetilde{w},\nu)_H=
(A\xi+A\nabla\widetilde{u}_\xi-\widetilde{A}^*\xi,\nu)_H\nonumber\\
&=&(A\nabla\widetilde{u}_\xi-\widetilde{A}^*\xi-A\nabla{u}_\xi+{A}^*\xi,\nu)_H=
(A\nabla\widetilde{u}_\xi-A\nabla{u}_\xi,\nu)_H+
\vert Y\vert\,\nu^T({A}^*-\widetilde{A}^*)\xi.\nonumber
\end{eqnarray}
Then
\begin{equation}\nonumber
 \vert (d_C,\nu)_H\vert \le 
 \lambda_2(A)\vert Y\vert^{\frac{1}{2}} \,\Vert \nabla\widetilde{u}_\xi-\nabla u_\xi
 \Vert_H \Vert \nu\Vert_2
+\delta_2\vert Y\vert\Vert \nu\Vert_2\Vert\xi\Vert_2
\le
\delta_1 \lambda_2(A)\vert Y\vert^{\frac{1}{2}}   \Vert \nu\Vert_2
+\delta_2\vert Y\vert\Vert \nu\Vert_2
\end{equation}
and thus $\Vert d_C\Vert_H\le \delta_1\lambda_2(A)\vert Y\vert^{\frac{1}{2}} 
+\delta_2\vert Y\vert$.
For the estimate of the distance $d_W$ between the $H$-orthogonal projection of $w_H\in H$ onto $W_h$ and  its approximation $\widetilde{w}\in W_{0,h}$ we obtain a standard inequality $\Vert d_W\Vert_H\le c h$, where 
the real constant $c$ does not depend on $h$; see e.g.~\cite{Brenner}.
Now we estimate $\Vert w_H-w_{\widetilde{\alpha}}\Vert_H$.
Note that from~\eqref{AEj} for any $\alpha\in\mathbb{R}^3$ the exact minimizer in~\eqref{min2} is
\begin{equation}\nonumber
w_{{\alpha}}=A(A^*)^{-1}\alpha+A\nabla u_{(A^*)^{-1}\alpha}-\alpha.
\end{equation}
Then setting $\alpha:=\widetilde{\alpha}=\widetilde{A}^*\xi$ we have
\begin{eqnarray}
\Vert w_H-w_{\widetilde{\alpha}}\Vert_H&=&
\Vert A\xi+A\nabla\widetilde{u}_\xi-\widetilde{A}^*\xi
-(
A(A^*)^{-1}\widetilde{\alpha}+A\nabla u_{(A^*)^{-1}\widetilde{\alpha}}-\widetilde{\alpha})\Vert_H\nonumber\\
&=&
\Vert A\xi+A\nabla\widetilde{u}_\xi
-(
A(A^*)^{-1}\widetilde{A}^*\xi+A\nabla u_{(A^*)^{-1}\widetilde{A}^*\xi})\Vert_H\nonumber\\
&\le&
\Vert A\xi-A(A^*)^{-1}\widetilde{A}^*\xi\Vert_H+
\Vert A(\nabla\widetilde{u}_\xi-
\nabla u_{(A^*)^{-1}\widetilde{A}^*\xi})\Vert_H\nonumber\\
&=&
\Vert A(A^*)^{-1}(A^*-\widetilde{A}^*)\xi\Vert_H+
\Vert A(\nabla\widetilde{u}_\xi-\nabla u_\xi+\nabla u_\xi-
\nabla u_{(A^*)^{-1}\widetilde{A}^*\xi})\Vert_H
\nonumber\\
&\le&
\frac{\lambda_2(A)}{\lambda_1(A^*)}\vert Y\vert\delta_2+
\lambda_2(A)\delta_1+
\Vert A(\nabla u_\xi-
\nabla u_{(A^*)^{-1}\widetilde{A}^*\xi})\Vert_H\nonumber\\
&\le&
\frac{\lambda_2(A)}{\lambda_1(A^*)}\vert Y\vert\delta_2+
\lambda_2(A)\delta_1+\lambda_2(A)^{\frac{1}{2}}
\Vert A^{1/2}(\nabla u_\xi-
\nabla u_{(A^*)^{-1}\widetilde{A}^*\xi})\Vert_H.\nonumber
\end{eqnarray}
Since the minimizers $u_\xi$ in~\eqref{min1} depend on $\xi$
linearly, then denoting $\nu=(A^*)^{-1}\widetilde{A}^*\xi$
we have 
\begin{equation}\nonumber
\Vert A^{1/2}(\nabla u_\xi-
\nabla u_{(A^*)^{-1}\widetilde{A}^*\xi})\Vert_H^2=
\int_Y (A(\nabla u_\xi-\nabla u_\nu),\nabla u_\xi-\nabla u_\nu)_{\mathbb{R}^3}
\, {\rm d}x=\vert Y\vert (\xi-\nu)^TA^*(\xi-\nu).
\end{equation}
We also have
\begin{equation}\nonumber
\xi-(A^*)^{-1}\widetilde{A}^*\xi=(A^*)^{-1}(A^*-\widetilde{A}^*)\xi.
\end{equation}
Then
\begin{equation}\nonumber
\Vert w_H-w_{\widetilde{\alpha}}\Vert_H\le 
\frac{\lambda_2(A)\vert Y\vert}{\lambda_1(A^*)}\delta_2+
\lambda_2(A)\delta_1+\frac{\lambda_2(A)^{\frac{1}{2}}
\vert Y\vert^{\frac{1}{2}}}{\lambda_1(A^*)^{\frac{1}{2}}}\delta_2,
\end{equation}
and finally,
\begin{equation}\nonumber
\Vert \widetilde{w}-w_{\widetilde{\alpha}}\Vert_H\le 
\Vert \widetilde{w}-w_{H}\Vert_H+
\Vert {w}_H-w_{\widetilde{\alpha}}\Vert_H\le c_1\delta_1+c_2\delta_2+c h,
\end{equation}
where $c_1,c_2>0$ depend on exteremal eigenvalues of $A$ and $A^*$,
and $c$ is provided by approximation properties of $W_{0,h}$ in $W_0$.
\end{proof}
}

\section{Numerical examples and implementation remarks}\label{sec5}

This section is devoted to the practical computation of the 
bounds on $A^*$. In the first part, we suggest a procedure 
for getting the upper and lower bounds on $A^*$.
In the second part, we go into detail and discuss
some implementation aspects and present some numerical examples.

\subsection{Algorithm}

The algorithm for providing two-sided bounds on $A^*$ is as follows:

{\bf Algorithm~1.} { Computing $A^*_h$ and $B^*_h$.}
\begin{enumerate}
\item[1.] Choose a FE discretization and solve~\eqref{min1} 
with $\xi=e^{(j)}$,  and get the solutions $u^{(j)}\in U$, $j=1,2,3$. \\
Set 
$(A_h^*)_{jk}=\frac{1}{\vert Y\vert}\int_Y(A(e^{(k)}+\nabla u^{(k)}),
e^{(j)}+\nabla u^{(j)})_{\mathbb{R}^3}\,\dd x$.
\item[2.] 
Choose a FE discretization and solve~\eqref{min2} 
with $\alpha=e^{(j)}$,  and get the solutions $w^{(j)}\in W_0$, $j=1,2,3$. \\
Set 
$(B_h^*)_{jk}=\frac{1}{\vert Y\vert}\int_Y(A^{-1}(e^{(k)}+w^{(k)}),
e^{(j)}+w^{(j)})_{\mathbb{R}^3}\,\dd x$.
\item[3.] The bounds on the true $A^*$ are $(B^*_h)^{-1}\preccurlyeq A^*
\preccurlyeq A^*_h$.
\end{enumerate}

\begin{rmrk}
\begin{enumerate}
\item[(a)]
 The choices of the FEM settings in steps 1 and 2 in
Algorithm~1 to get $B_h^*$ and $A_h^*$, respectively, are independent.
\item[(b)]
 Matrices $(B^*_h)^{-1}$ and $A_h^*$ are the lower and upper bounds,
respectively, to $A^*$, only if the system matrices and the right-hand side vectors are evaluated exactly. In practice,
this means that the elements of the system matrices must be integrated 
exactly; otherwise, the bounds are not guaranteed anymore.
On the other hand, the solutions $u^{(j)}$ and $w^{(j)}$ do not
have to be exact solutions of the discretized systems; they only have to belong to the respective function spaces,
i.e.~$u^{(j)}\in U$ and $w^{(j)}\in W_0$, $j=1,2,3$.
\item[(c)] Taking any basis $\{\xi^{(1)},\xi^{(2)},\xi^{(3)}\}$
of $\mathbb{R}^3$ and defining $A^*_h$ using any functions 
$u^{(1)},u^{(2)},u^{(3)}\in U$ in~\eqref{min1}
yields $A^*{\preccurlyeq}  A^*_h$.
The closest upper bound is, however, obtained by such $u^{(j)}$,
$j=1,2,3$, that minimize~\eqref{min1}.
An analogous statement holds for the lower bounds.
\end{enumerate}
\end{rmrk}

In practical computation, instead of solving~\eqref{min2} 
for { a given } $\alpha\in \mathbb{R}^3$
in some subspace of $W_0$, 
we can use formula~\eqref{AEj} in {Theorem}~\ref{lem55} to get
a good approximation 
$\widetilde w$ of its minimizing field ${w_\alpha}\in W_0$ and thus to get an approximation {$\widetilde{B}_h^*$ of ${B}^*$},
\begin{equation}\label{w01}
({\widetilde{B}_h^*}\alpha,\alpha)_{{\mathbb R}^3}=
\frac{1}{\vert Y\vert}\int_Y(  A^{-1}(\alpha+\widetilde w),
\alpha+\widetilde w)_{{\mathbb R}^3}\, \dd x.
\end{equation}
{ In practice, we cannot use 
$w_1=A(\xi+\nabla u_\xi)-A^*\xi$ as $\widetilde{w}$
directly, because we know neither $u_\xi$ nor $A^*$
exactly. Instead, we must only consider some approximations
of these quantities $\widetilde{u}_\xi$ and  $\widetilde{A}^*$,
respectively, obtained from step~1 of Algorithm~1,
and then set
\begin{equation}
    w_H=A(\xi+\nabla \widetilde{u}_\xi)-\widetilde{A}^*\xi\in H.
\end{equation}
Such a $w_H$, however, may not belong to $W_0$ in general. 
 Therefore, we compute a projection $\widetilde{w}$ of
 $w_H$ onto some finite-dimensional space $W_{0,h}\subset W_0$.
  In this paper, we propose to use the orthogonal 
  projection on a FE space $W_{0,h}\subset W_0$.
The orthogonality is understood in the sense of 
$H=L^2(Y,\mathbb{R}^3)$ inner product.
Then $\widetilde{w}\in W_{0,h}$ is obtained as a solution of
\begin{equation}\label{P1}
    (\widetilde{w},v)_{H}=
    ( w_H,v  )_{H},\quad v\in W_{0,h}.
\end{equation}
and $\widetilde{B}_h^*$ is then obtained from~\eqref{w01}
where we set $\alpha:=\widetilde{\alpha}=\widetilde{A}^*\xi$.
Since $\widetilde{w}$  is not the minimizer of~\eqref{min2}
with $\widetilde{\alpha}=\widetilde{A}^*\xi$, we have 
\begin{equation}
    ({B}^*_h\widetilde{\alpha},\widetilde{\alpha})_{\mathbb{R}^3}\le
    (\widetilde{B}^*_h\widetilde{\alpha},\widetilde{\alpha})_{\mathbb{R}^3}.
\end{equation}  
In section~\ref{secnumsol} we prove that the approximations $A^*_h$ and $(\widetilde{B}_h^*)^{-1}$
of the exact matrix $A^*$ as its upper and lower bounds, respectively, are asymptotically of the same rate.
  We can also consider any other kind of projection, such as
  energy-like inner product projections. In our numerical experiments, the simplest one~\eqref{P1} gives 
  good numerical results; thus, we do not involve the others in
  our consideration.
  It is important to emphasize that due to periodic boundary conditions and using
uniform grids, the linear system matrix is block circulant,
and thus the solution of~\eqref{P1} can be obtained by 
the fast discrete Fourier transform without solving any 
system of linear equations.  \\

  {\bf Algorithm~2.} Computing $\widetilde{B}^*_h$ from approximate solutions of the primal problem.
  \begin{enumerate}
\item[1.] 
       Using Algorithm~1, step 1, for       
      $\xi=e^{(j)}\in\mathbb{R}^3$ get (possibly inexact) solutions $\widetilde{u}^{(j)}\in U$, $j=1,2,3$, and 
      an approximation $\widetilde{A}^*$ to $A^*\in\mathbb{R}^{3\times 3}$.
\item[2.] 
For $j=1,2,3$:\\
set $w^{(j)}_H=A\left(e^{(j)}+\nabla \widetilde{u}^{(j)}\right)-\widetilde{A}^*e^{(j)}$\\
find $\widetilde{w}^{(j)}\in W_{0,h}$ such that $(\widetilde{w}^{(j)},v)_H=
(w^{(j)}_H,v)_H$, for all $v\in W_{0,h}$\\
set $\alpha^{(j)}=\widetilde{A}^*e^{(j)}$
\item[3.]
Build $\widetilde{B}_h^*$ such that
$(\widetilde{B}_h^*\alpha^{(k)},\alpha^{(j)})_{\mathbb{R}^3}=
\frac{1}{\vert Y\vert}\int_Y \left(A^{-1}(\alpha^{(k)}+\widetilde{w}^{(k)}),\alpha^{(j)}+\widetilde{w}^{(j)}
\right)_{\mathbb{R}^3}\,\dd x$.
  \end{enumerate}
Due to the regular grid of $Y$ and periodic boundary conditions,
the mass matrix of the linear system arising in~\eqref{P1}
is block circulant, which enables using the discrete fast Fourier transform
to obtain the solution without solving any system of linear equations.
}

\subsection{Numerical examples}\label{numex}

In our numerical experiments, we choose $a_i=2\pi$, $i=1,2,3$
and consider $N_{\rm vox}=N_1N_2N_3$ voxels in $Y$.
Each voxel is split into six tetrahedral elements of the same volume,
thus the number of elements is $N_{\rm ele}=6N_{\rm vox}$. 
We use continuous piece-wise linear {Lagrange} FE basis functions.
To better pronounce the properties of the algebraic objects involved
in our computation and/or to obtain insight into a matrix-free
{form} of the algorithms, we can (in the spirit of e.g.~\cite{LadeckyEtAl22})
construct the stiffness matrices and the right-hand side vectors from 
the so-called derivative matrices and from the coefficients in the following way.
We build derivative matrices $\mathsf D_{j}$, $j=1,2,3$, which map the vector of nodal values of the {basis} function{~$\phi_m$} to the values of derivatives in quadrature points $x^q$~(centers of elements),
$(\mathsf{D}_j)_{km}=\frac{\partial \phi_m}{\partial x_j}(x^q_k)$,
$k=1,\dots,N_{\rm ele}$, $m=1,\dots,N_{\rm vox}$, $j=1,2,3$.
Following~\eqref{Q123}, we denote
\begin{equation}\label{D}
 {\mathsf D}_{\rm grad}=\left(\begin{array}{c}{\mathsf D}_1\\
    {\mathsf D}_2\\ {\mathsf D}_3\end{array}\right),\quad
    {\mathsf D}_{{\rm curl}}=\left(\begin{array}{ccc}
    0&{\mathsf D}_3&-{\mathsf D}_2\\
    -{\mathsf D}_3&{\mathsf 0}&{\mathsf D}_1\\ 
    {\mathsf D}_2 & -{\mathsf D}_1&0\end{array}\right),
\end{equation}\nonumber
{i.e., ${\mathsf D}_{{\rm curl}}=(Q_3 {\mathsf D}_{\rm grad},Q_2 {\mathsf D}_{\rm grad},Q_1 {\mathsf D}_{\rm grad})$.
}
The stiffness matrices and the right-hand sides of the discretized primal and dual problems~\eqref{min1} and~\eqref{min2}, respectively, are
\begin{equation*}
{\mathsf K}={\mathsf D}_{\rm grad}^T\mathsf{AD}_{\rm grad},\;
 \mathsf{b}=-\mathsf{D}_{\rm grad}^T\mathsf{Ae}^\xi, \qquad
{\mathsf K}_{{\rm dual}}={\mathsf D}_{{\rm curl}}^T\mathsf{A}_{\rm inv}
\mathsf{D}_{{\rm curl}},\; \mathsf{b}_{{\rm dual}}=-\mathsf{D}_{{\rm curl}}^T\mathsf{A}_{\rm inv}\mathsf{e}^\alpha,
\end{equation*}
where  
$\mathsf{A}\in\mathbb{R}^{3N_{\rm ele}\times 3N_{\rm ele}}$ is composed from $3\times 3$ blocks 
which are $N_{\rm ele}\times N_{\rm ele}$ diagonal matrices
with the elements of $A(x_k^q)\in \mathbb{R}^{3\times 3}$
as $k$-th elements of all respective nine diagonals.
Similarly, $\mathsf{A}_{\rm inv}\in\mathbb{R}^{3N_{\rm ele}\times 3N_{\rm ele}}$ involves the elements 
of $(A(x_k^q))^{-1}\in \mathbb{R}^{3\times 3}$.
Vectors $\mathsf{e}^\xi,\mathsf{e}^\alpha\in\mathbb{R}^{3N_{\rm ele}}$
contain repeating elements of $\xi$ or $\alpha$, respectively.
Specifically, $\mathsf{e}^\xi={\rm kron}\,(\xi,\mathsf{1})$
or $\mathsf{e}^\alpha={\rm kron}\,(\alpha,\mathsf{1})$ when
we consider $\xi$ in~\eqref{min1} or $\alpha$ in~\eqref{min2}, respectively,
where $\rm kron$ denotes the Kronecker product and $\mathsf{1}=(1,\dots,1)^T
\in \mathbb{R}^{N_{\rm ele}}$.

After obtaining the solutions $\mathsf u$ of $\mathsf{K}\mathsf{u}=\mathsf{b}$
and $\mathsf w$ of 
{ $\mathsf{K}_{{\rm dual}}\mathsf{w}=\mathsf{b}_{{\rm dual}}$, the approximate homogenized matrices follow from}
\begin{eqnarray}
\xi^TA^*_h\xi &=& ( {\mathsf D}_{\rm grad}\mathsf{u}+\mathsf{e}^\xi)^T\mathsf{A}( {\mathsf D}_{\rm grad}\mathsf{u}+\mathsf{e}^\xi)
\label{Ah1}\\
\alpha^TB^*_{h}\alpha &=& 
({\mathsf D}_{{\rm curl}}\mathsf{w}+\mathsf{e}^\alpha)^T
\mathsf{A}_{\rm inv}({\mathsf D}_{{\rm curl}}\mathsf{w}+\mathsf{e}^\alpha),
\label{Bh1}
\end{eqnarray}
{for arbitrary $\xi$ and $\alpha$ in $\mathbb{R}^3$.}

\begin{rmrk}[Dual systems]\label{rmk:dual_system}
 The dimension of the stiffness {matrix $\mathsf{K}_{{\rm dual}}$  for getting the lower bound
is $3N_{\rm vox}$ in contrast to the dimension $N_{\rm vox}$
of the stiffness matrix $\mathsf{K}$ of the primal problem (for getting the upper bound)}.
\end{rmrk}

For both experiments, we set $N_1=N_2=N_3$ and consider the positive definite matrices $A$ that are {piecewise constant functions. Furthermore, the discretization is chosen such that $A$ 
is constant on each element.} {
Then the numerical quadrature rule with one quadrature point per element yields that the stiffness matrices and the right-hand side vectors are computed exactly.}
The tetrahedra in the mesh are oriented such that only the vertices $(0,0,0)$ and $(1,1,1)$ are connected by edges of the tetrahedra with seven other vertices, while the remaining vertices are connected each with four other vertices. The CG algorithm was stopped when the residual of the initial guess $\mathsf{x} = \mathsf{0}$ was reduced by the factor $10^{-9}$. {All results can be reproduced using codes available in a~Zenodo repository~\cite{Gaynutdinova:2023:zenodo}}.

{
 \paragraph{Example~1} For the first example, we consider the coefficients in the form}
\begin{eqnarray}
   A(x)=\left(\begin{array}{ccc}
   7+{\rm sign}(s_1s_2)  &  -2-{\rm sign}(s_2s_3)  &  {\rm sign}(s_1s_2s_3)  \\
   -2- {\rm sign}(s_2s_3)   & 4.01+{\rm sign}(s_1s_2)  &   0\\
    {\rm sign}(s_1s_2s_3)  & 0 & 3+{\rm sign}(s_2s_3)
   \end{array}\right), \quad
   x \in Y, \quad
   s_j=\sin(\textstyle\frac{3}{2}x_j).
   \label{eq:example1}
\end{eqnarray}

{The data collected in Table~\ref{tab:Example1} validate the upper-lower structure of the bounds, $(\widetilde{B}_{h}^*)^{-1}\preccurlyeq(B_{h}^*)^{-1}\preccurlyeq A_h^*$, as follows from the positive definiteness of the differences $A_h^* - (B_{h}^*)^{-1}$ and $\widetilde{B}_{h}^* - B_{h}^*$, and the monotone convergence of the diagonal elements upon the mesh refinement. We notice that for the primal problem, the number of CG iterations doubles with doubling the number of voxels per cell edge. Because of the expanded solution space, recall Remark~\ref{rmk:dual_system}, the dual problem needs substantially more iterations to converge as the number of iterations increases four times under the mesh refinement. The lower bound determined from the projected solution $(\widetilde{B}_{h}^*)^{-1}$, on the other hand, provides an excellent approximation at a negligible computational cost. Additionally, we observe that the accuracy of the bounds surpasses the values reported for the Fourier basis in~\cite{VondrejcEtAl2015,Vondrejc2016}; even rather coarse discretization of $N_1 = N_2 = N_3 = 24$ yields the relative difference between upper and lower bounds less than $1\%$; see also~\cite{Vondrejc2020} for a related study for the primal problem. This result is encouraging for the target applications in computational micromechanics of heterogeneous materials.} 

\begin{table}[h]
$N_1 = N_2=N_3=6$ \\
\begin{tabular}{ccc}
${(\widetilde{B}_{h}^*)^{-1}}$ &&\\
\hline
  6.5702  & -2.1432 &  -0.0629  \\  
   -2.1432  &  3.8983 &  -0.0096 \\  
   -0.0629  & -0.0096 &   2.7496 \\  
\hline
   &     &    \\
\end{tabular}\hskip10pt
\begin{tabular}{ccc}
$(B_{h}^*)^{-1}$ &&\\
\hline
 6.6193  & -2.1350  & -0.0562\\
   -2.1350  &  3.9140 &  -0.0064\\
   -0.0562 &  -0.0064  &  2.7756\\   
\hline
90  &  89  &  89\\
\end{tabular}\hskip10pt
\begin{tabular}{ccc}
$A_h^*$ &&\\
\hline
  6.9126 &  -2.0937  & -0.0114\\
   -2.0937 &   4.0453 &  -0.0029\\
   -0.0114  & -0.0029  &  2.9602\\  
\hline
 36   & 35 &   36\\
\end{tabular}   

\textcolor{red}{
$
\mathrm{eig}\, ( A_h^* - ( B_{h}^*)^{-1} )
=
(0.1205, 0.1707, 0.3181),
\hspace{1ex}
\mathrm{eig}\,  ( ( B_{h}^*)^{-1} - ( \widetilde{B}_{h}^*)^{-1} )
=
(0.0135, 0.0243, 0.0528) 
$}

\medskip $N_1 = N_2=N_3=12$

\begin{tabular}{ccc}
$(\widetilde{B}_{h}^*)^{-1}$ &&\\
\hline
6.7067  & -2.1203  & -0.0471    \\
   -2.1203  &  3.9621  & -0.0083 \\  
   -0.0471 &  -0.0083  &  2.8249 \\  
\hline
  &    &    \\   
\end{tabular}\hskip10pt
\begin{tabular}{ccc}
$(B_{h}^*)^{-1}$ &&\\
\hline 
 6.7239  & -2.1171 & -0.0437\\
   -2.1171  &  3.9675 &  -0.0073\\
   -0.0437  & -0.0073  &  2.8367  \\
\hline
361 &  360 &  364\\
\end{tabular}\hskip10pt
\begin{tabular}{ccc}
$A_h^*$ &&\\
\hline   
     6.8414  & -2.1012 &  -0.0253\\
   -2.1012 &   4.0189 &  -0.0051\\
   -0.0253  & -0.0051  &  2.9105\\    
\hline
74  &  73  &  75 \\
\end{tabular}

\textcolor{red}{
$
\mathrm{eig}\, ( A_h^* - ( B_{h}^*)^{-1} )
=
(0.0475, 0.0677, 0.1275),
\hspace{1ex}
\mathrm{eig}\,  ( ( B_{h}^*)^{-1} - ( \widetilde{B}_{h}^*)^{-1} )
=
(0.0047, 0.0102, 0.0197) 
$}

\medskip $N_1 = N_2 = N_3 = 24$

\begin{tabular}{ccc}
$(\widetilde{B}_{h}^*)^{-1}$ &&\\
\hline
6.7625& -2.1117 & -0.0390\\
-2.1117& 3.9867& -0.0073\\
-0.0390& -0.0073&  2.8594\\
\hline
     &          &       \\   
\end{tabular}\hskip10pt
\begin{tabular}{ccc}
$(B_{h}^*)^{-1}$ &&\\
\hline 
 6.7683   & -2.1106 &-0.0378\\
  -2.1106  &3.9885&-0.0070\\
 -0.0378 &-0.0070&2.8636\\ 
\hline
1,404    &    1,398   &     1,410 \\
\end{tabular}\hskip10pt
\begin{tabular}{ccc}
$A_h^*$ &&\\
\hline   
 6.8091  &-2.1049&-0.0314\\
 -2.1049 & 4.0063& -0.0060 \\
 -0.0314 & -0.0060 &2.8891\\
\hline
158     &   156    &    157 \\
\end{tabular}

\textcolor{red}{
$
 \mathrm{eig}\, ( A_h^* - ( B_{h}^*)^{-1} )
=
(0.0164, 0.0234, 0.0444),
\hspace{1ex}
\mathrm{eig}\,  ( ( B_{h}^*)^{-1} - ( \widetilde{B}_{h}^*)^{-1} )
=
(0.0015, 0.0036, 0.0066) 
$}

\bigskip
\caption{
{Bounds on and estimates of the diagonal and off-diagonal elements, respectively, of the homogenized matrix $A^*$ for coefficients specified in Example~1, Eq.~\eqref{eq:example1}. $A_h^*$ and $(B_{h}^*)^{-1}$ denote the upper and lower bounds determined according to relations~\eqref{Ah1} and~\eqref{Bh1}, and $(\widetilde{B}_{h}^*)^{-1}$ is the lower bound determined from the projected solution~\eqref{P1}. $N_j$, for $j = 1, 2, 3$, stands for the number of voxels per edge length, and the bottom rows collect the number of iterations of the conjugate gradient (CG) method, noting that computing $(\widetilde{B}_{h}^*)^{-1}$ does not involve any CG iterations.}}
\label{tab:Example1}
\end{table}

\begin{figure}[h]
   \centering
   \includegraphics[width=\textwidth]{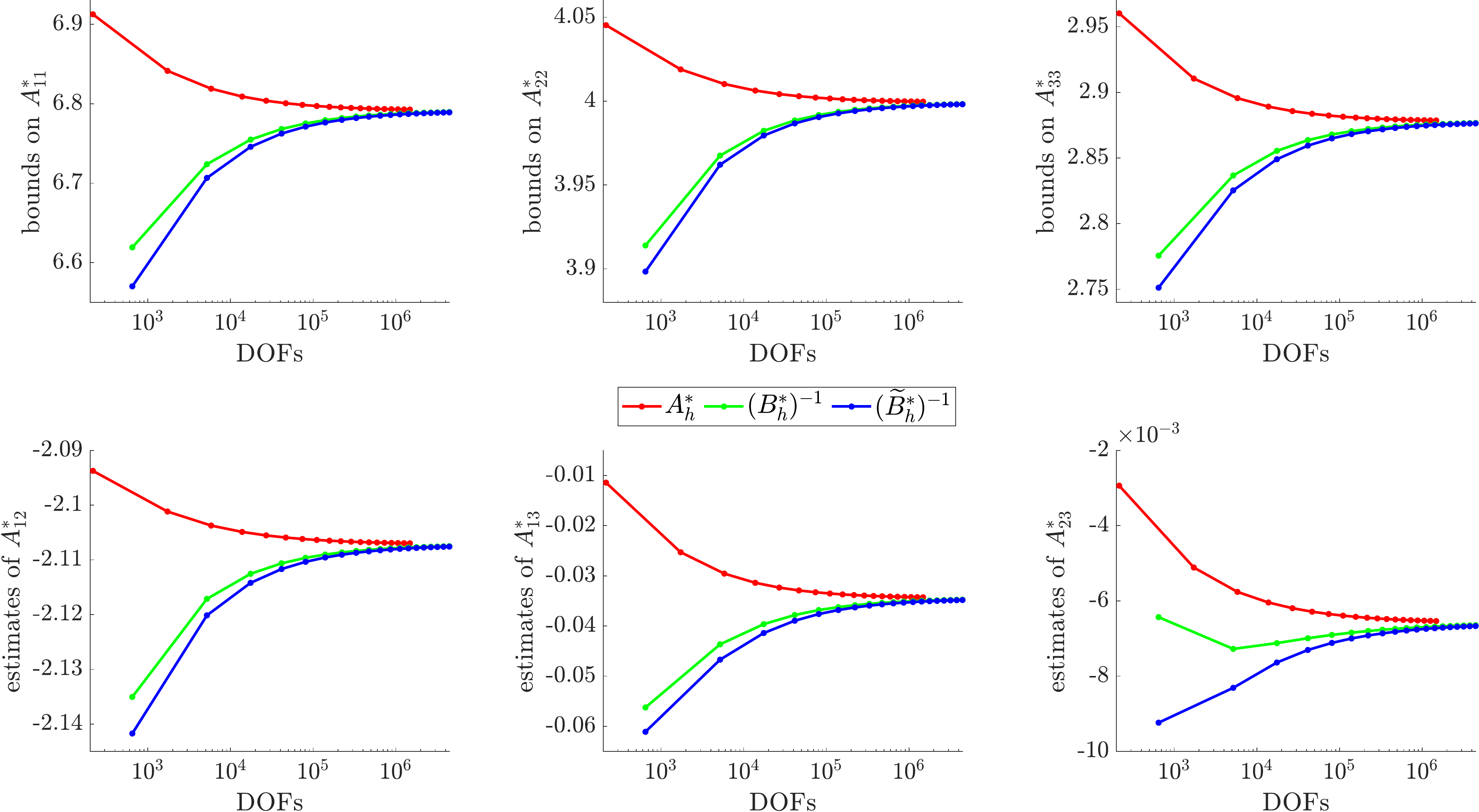}
   \caption{{Estimates of elements of $A^*$ for Example~1, Eq.~\eqref{eq:example1}.
   First row: Guaranteed upper and lower bounds on diagonal elements of $A^*$ obtained as 
   the respective diagonal elements of     
   $A^*_h$, $B^*_h$ and $\widetilde{B}^*_h$. 
   Second row: Analogously obtained estimates (not guaranteed bounds) for off-diagonal elements $A^*_{12}$, $A^*_{13}$, and $A^*_{23}$. Matrices $A^*_h$ and $B^*_h$ are computed
   according to Algorithm~1; matrices $\widetilde{B}^*_h$
   are computed according to Algorithm~2.
   The horizontal axes denote the total number of degrees of freedom (DOFs) for every problem.}}
   \label{fig1}
\end{figure}

{
Figure~\ref{fig1} provides further insights into the behavior of the elements of upper-lower bounds matrices with an increasing number of degrees of freedom. Indeed, these extended data confirm the convergence of all the elements preserving the ordering $(\widetilde{B}_{h}^*)^{-1}\preccurlyeq(B_{h}^*)^{-1}\preccurlyeq A_h^*$, with analogous ordering for all diagonal elements as a consequence. Figure~\ref{fig1}  additionally certifies the quality of the projection-based lower bound $(\widetilde{B}_{h}^*)^{-1}$ and the extra computational cost of the lower bound $(B_{h}^*)^{-1}$ that is reflected in its horizontal shift relative to the upper bound $A_h^*$.}

\begin{figure}[h]
   \centering
   \includegraphics[width=\textwidth]{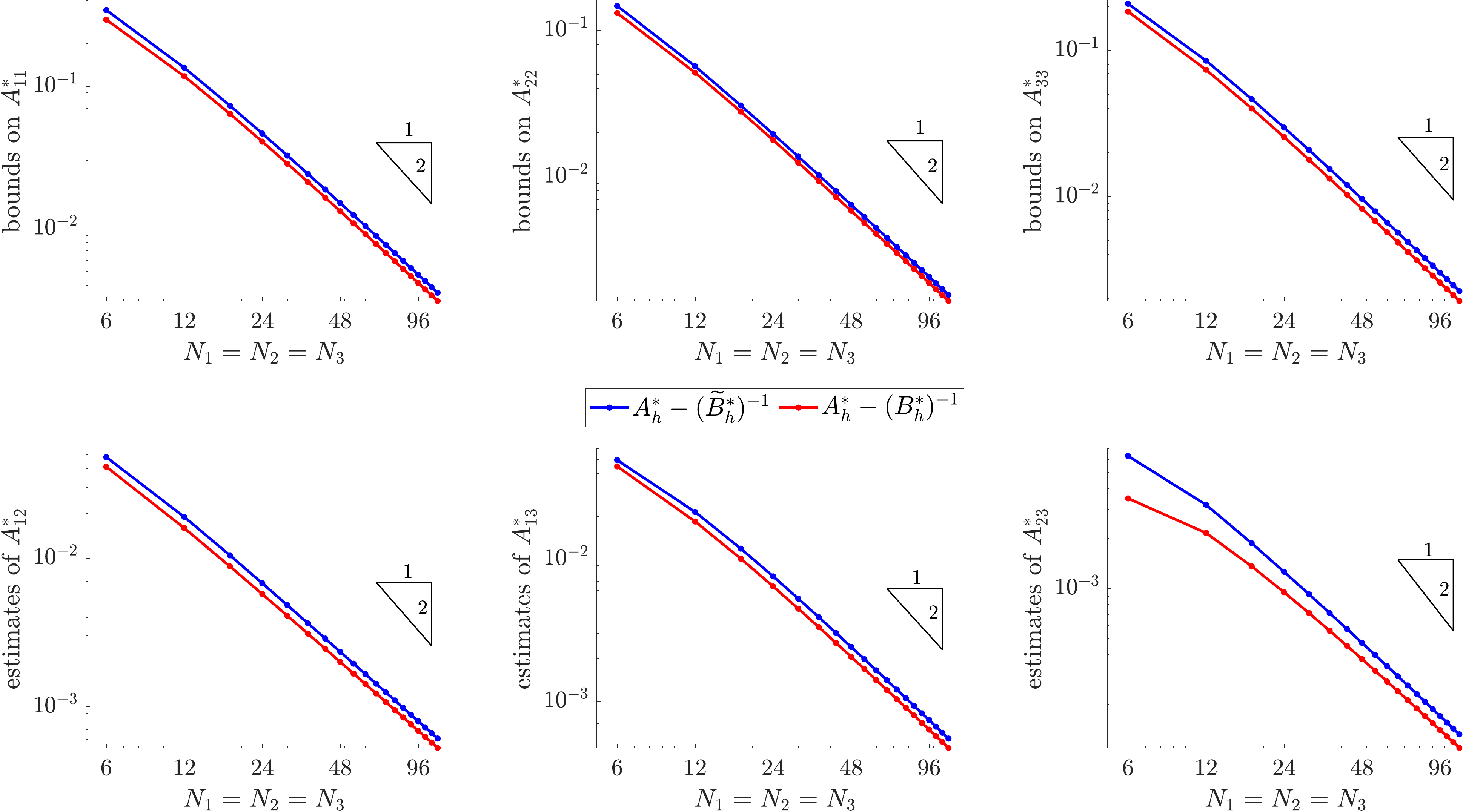}   
   \caption{Convergence rates for the error estimates of the homogenized properties for coefficients of Example~1, Eq.~\eqref{eq:example1}. First row: Guaranteed error bounds on diagonal elements; second row: Estimates (not guaranteed bounds) on off-diagonal elements. $N_j$, $j = 1, 2, 3$ stands for the number of voxels per cell $Y$ edge, and the triangle indicates the quadratic convergence rate with the size $h$, i.e., $h^2 \approx N_1^{-2}$.}
   \label{fig4}
\end{figure}

{
Finally, we examine the rate of convergence of the homogenized properties by plotting the elements of the error estimates $( A_h^* - (B_{h}^*)^{-1})$ and $( A_h^* - (\widetilde{B}_{h}^*)^{-1} )$ against the number of voxels per edge $N_1 = N_2 = N_3$ in the double-logarithmic scale in Figure~\ref{fig4}. The results confirm the asymptotic scaling of $N_1^{-2} \approx h^2$, consistently with the results of Lemma~\ref{lem_dist} and the approximation property of piecewise linear basis functions, e.g.,~\cite{Brenner} (for the estimate $( A_h^* - (B_{h}^*)^{-1})$) and the additional results of Lemma~\ref{lem:projection_error} (for $( A_h^* - (\widetilde{B}_{h}^*)^{-1} )$). Note that the quadratic rate of convergence corresponds to the fact that the discontinuities in the coefficient $A(x)$ from Eq.~\eqref{eq:example1} are exactly captured with the finite element mesh. The general case of discontinuity-unfitted meshes, typical of voxel-based geometries, corresponds to the linear convergence rate in homogenized properties; see~\cite{Schneider2022nme,Ye:2023} for further discussion.}

{
\paragraph{Example~2} To verify that the conclusions drawn from Example 1 also hold for other types of problem data, we perform an analogous analysis for the following coefficients:}
\begin{eqnarray}\label{eq:example2}
   A(x)=(2+{\rm sign}(s_1s_2s_3))I, \quad 
   x=(x_1,x_2,x_3)^T\in Y, \quad 
   s_j=\sin(\textstyle\frac{3}{2}x_j)    
\end{eqnarray}
{
The results on the upper-lower bounds for discretizations $N_1 = N_2 = N_3 = 6, 12$, and $24$, collected in Table~\ref{tab:Example2}, confirm identical behavior as reported for Example 1 in Table~\ref{tab:Example1}. Similarly, the convergence graphs exhibit a~similar pattern, although they are omitted here for brevity. Notice that the isotropy of local coefficients and the symmetry of their spatial distribution visible from Eq.~\eqref{eq:example2} imply that the off-diagonal elements of the effective matrix $A^*$ are zero. In contrast, the corresponding values of both matrices $A^*_h$ and $B^*_h$ are non-positive. Therefore, they estimate the exact value rather than bounding it.}

\begin{table}[h]

$ N_1 = N_2 = N_3 = 6$

\begin{tabular}{ccc}
$(\widetilde{B}_{h}^*)^{-1}$ &&\\
\hline
 1.7035 &  -0.0043  & -0.0043  \\  
   -0.0043  &  1.7035 &  -0.0043 \\  
   -0.0043 &  -0.0043  &  1.7035 \\    
\hline
   &      &  \\
\end{tabular}\hskip10pt
\begin{tabular}{ccc}
$(B_{h}^*)^{-1}$ &&\\
\hline
  1.7066 &  -0.0043 &  -0.0043 \\   
   -0.0043  &  1.7066 &  -0.0043  \\ 
   -0.0043 &  -0.0043 &   1.7066  \\
\hline
 65   & 65 &   65 \\
\end{tabular}\hskip10pt
\begin{tabular}{ccc}
$A_h^*$ &&\\
\hline
 1.9446  & -0.0016 &  -0.0016\\
 -0.0016 &   1.9446 &  -0.0016\\
 -0.0016 &  -0.0016 &   1.9446\\
\hline
 28  &  28  &  28\\
\end{tabular}

\textcolor{red}{
$
\mathrm{eig}\, ( A_h^* - ( B_{h}^*)^{-1} )
=
(0.2353, 0.2353, 0.2434),
\hspace{1ex}
\mathrm{eig}\,  ( ( B_{h}^*)^{-1} - ( \widetilde{B}_{h}^*)^{-1} )
=
(0.0030, 0.0031, 0.0031) 
$}

\medskip $N_1 = N_2 = N_3 = 12$

\begin{tabular}{ccc}
$(\widetilde{B}_{h}^*)^{-1}$ &&\\
\hline
   1.7831  & -0.0023 &  -0.0023 \\  
   -0.0023  &  1.7831 &  -0.0023 \\  
   -0.0023  & -0.0023 &   1.7831  \\ 
\hline
  &    &    \\
\end{tabular}\hskip10pt
\begin{tabular}{ccc}
$(B_{h}^*)^{-1}$ &&\\
\hline
   1.7859 &  -0.0022 &  -0.0022\\
   -0.0022  &  1.7859 &  -0.0022\\
   -0.0022 &  -0.0022 &   1.7859\\
\hline
 252 & 252 &  252\\
\end{tabular}\hskip10pt
\begin{tabular}{ccc}
$A_h^*$ &&\\
\hline
  1.8938 &  -0.0002  & -0.0002\\
   -0.0002  &  1.8938 &  -0.0002\\
   -0.0002 &  -0.0002   & 1.8938\\
\hline
65 &   65 & 65\\
\end{tabular}

\textcolor{red}{
$
\mathrm{eig}\, ( A_h^* - ( B_{h}^*)^{-1} )
=
(0.1059, 0.1059, 0.1119),
\hspace{1ex}
\mathrm{eig}\,  ( ( B_{h}^*)^{-1} - ( \widetilde{B}_{h}^*)^{-1} )
=
(0.0028, 0.0028, 0.0029)$ }

\medskip $N_1 = N_2 = N_3 = 24$

\begin{tabular}{ccc}\label{table:Examples}
$(\widetilde{B}_{h}^*)^{-1}$ &&\\
\hline
   1.8214 &  -0.0008 &  -0.0008  \\  
   -0.0008  &  1.8214 &  -0.0008 \\  
   -0.0008 &  -0.0008 &   1.8214 \\ 
\hline
    &   &   \\
\end{tabular}\hskip10pt
\begin{tabular}{ccc}
$(B_{h}^*)^{-1}$ &&\\
\hline
   1.8231  & -0.0008 &  -0.0008 \\   
   -0.0008  &  1.8231  & -0.0008 \\  
   -0.0008 &  -0.0008 &   1.8231 \\  
\hline
974  &974 &974\\
\end{tabular}\hskip10pt
\begin{tabular}{ccc}
$A_h^*$ &&\\
\hline
1.8671 & -0.0000 & -0.0000\\
-0.0000  &  1.8671 &  -0.0000\\
-0.0000 &  -0.0000  &  1.8671\\
\hline
131 & 131   & 131\\
\end{tabular}

\textcolor{red}{
$
\mathrm{eig}\, ( A_h^* - ( B_{h}^*)^{-1} )
=
(0.0433, 0.0433, 0.0456),
\hspace{1ex}
\mathrm{eig}\,  ( ( B_{h}^*)^{-1} - ( \widetilde{B}_{h}^*)^{-1} )
=
(0.0017, 0.0017, 0.0017) 
$}

\bigskip

\caption{
{Bounds on and estimates of the diagonal and off-diagonal elements, respectively, of the homogenized matrix $A^*$ for coefficients specified in Example~2, Eq.~\eqref{eq:example2}. $A_h^*$ and $(B_{h}^*)^{-1}$ denote the upper and lower bounds determined according to relations~\eqref{Ah1} and~\eqref{Bh1}, and $(\widetilde{B}_{h}^*)^{-1}$ is the lower bound determined from the projected solution~\eqref{P1}. $N_j$, for $j = 1, 2, 3$, stands for the number of voxels per edge length, and the bottom rows collect the number of iterations of the conjugate gradient (CG) method, noting that computing $(\widetilde{B}_{h}^*)^{-1}$ does not involve any CG iterations.   
}}
\label{tab:Example2}
\end{table}

\section{Conclusions}\label{sec6}

We have presented a reliable numerical procedure for computing the guaranteed
upper and lower bounds on a~homogenized coefficient matrix defined by the cell problem~\eqref{min1}. While the upper bound naturally follows from solving~\eqref{min1} approximately, e.g.~on a finite-dimensional subspace of $H^1_{\rm per}(Y,\mathbb{R})$, the lower bound requires building conforming approximations of the dual problem defined on $H^{{\rm div},0}_{\rm per}(Y,\mathbb{R}^3)$.

Let us note that an efficient numerical method can be rarely developed without employing all relevant {
sources} of the problem, starting from the physical meaning and formulation, through the choice of approximation spaces
up to the numerical solution method and appropriate preconditioning; see a more involved exposition in~\cite{Malek}. Therefore, we presented most of this path here within our suggestion on the estimation method by proving the dual problem (Lemma~\ref{Lem1} and~\ref{lem55}) based on a simple optimization result (Lemma~\ref{lem2}). We showed the Helmholtz decomposition of periodic
fields and introduced another representation of $H^{{\rm div},0}_{\rm per}(Y,\mathbb{R}^3)$ (Lemma~\ref{Lem11}) which led to approximation spaces based on FE discretization. The lower bounds were then obtained by the numerical solution of~\eqref{min2}.

{
Since the dual problem~\eqref{min2} is more computationally demanding than the primal one, it 
can be worth finding some good approximation of the minimizer of~\eqref{min2}. We show how such a minimizer can
be obtained using just a projection evaluated by
the discrete fast Fourier transform.

Our next research will focus on 
efficient preconditioner for the problem~\eqref{min2}, building on recent developments in the Laplace preconditioning of the primal problem~\cite{Gergelits2019,Gergelits2020,Ladecky2021,Pultarova2021}, to mitigate the dependence of the number of iterations to convergence on spatial discretization (recall Tables~\ref{tab:Example1} and~\ref{tab:Example2}), and} on obtaining the two-sided bounds on more involved problems in elasticity. 

\begin{acknowledgement}
{The authors thank 
Michal K\v r\' i\v zek~(Institute of Mathematics, Czech Academy of Sciences) for valuable discussions
and both anonymous referees for 
insightful and helpful comments.}
All authors acknowledge funding by the European Regional Development Fund (Centre of Advanced Applied Sciences -- CAAS, CZ 02.1.01/0.0/0.0/16\_019/0000778), by  the Czech Science Foundation (projects No.~20-14736S (ML), { and project No. 22-35755K (LG))}, 
the Student Grant Competition of the Czech Technical University in Prague {(project No.~SGS23/002/OHK1/1T/11 (LG)).}
\end{acknowledgement}

\ifx\undefined\bysame
\newcommand{\bysame}{\leavevmode\hbox to3em{\hrulefill}\,}
\fi

\appendix

\section{Proof of Theorem~\ref{charU}}\label{apxA}

Let us first introduce some auxiliary propositions.
Given $M\subset \mathbb{R}^3$, $c\in\mathbb{R}^3$ and $b\in\mathbb{R}$, let us use the notation $c+b M=\{c+bx;\;x\in M\}$ and recall that $Y=(0,a_1)\times (0,a_2) \times (0,a_3)$.
\begin{lmm}\label{elvhfbvhk1}
Let $u\in L^2(Y,{\mathbb R}^3)$. Then for all $Y$-periodic 
$\phi\in\mathcal{C}^{\infty}(\mathbb{R}^3,{\mathbb R})$ and
$c\in{\mathbb R}^3$ we have
\begin{align*}
&\int_{Y}u\cdot\nabla\phi\,\dd x=\int_{c+Y}u_{\rm per}\cdot\nabla\phi\,\dd x.
\end{align*}
\end{lmm}
\begin{proof}
Let us take first $c=(c_1,0,0)$ where $c_1=ka_1+b_1$, $k\in\mathbb{Z}$ and $\vert b_1\vert<a_1$. Denoting $b=(b_1,0,0)$ we have 
\begin{align*}
&\int_{c+Y}u_{\rm per}\cdot\nabla\phi\,\dd x=
\int_{b+Y}u_{\rm per}\cdot\nabla\phi\,\dd x
=\int_{(b+Y)\cap Y}u\cdot\nabla\phi\,\dd x+\int_{(b+Y)\setminus Y}u_{\rm per}\cdot\nabla\phi\,\dd x\\
&=\int_{Y\cap(b+Y)}u\cdot\nabla\phi\,\dd x+\int_{Y\setminus(b+Y)}u\cdot\nabla\phi\,\dd x=\int_{Y}u\cdot\nabla\phi\,\dd x.
\end{align*}
Then for $c=(c_1,c_2,c_3)\in \mathbb{R}^3$
\begin{align*}
&\int_{c+Y}u_{\rm per}\cdot\nabla\phi\,\dd x=
\int_{(0,c_2,c_3)+Y}u_{\rm per}\cdot\nabla\phi\,\dd x=
\int_{(0,0,c_3)+Y}u_{\rm per}\cdot\nabla\phi\,\dd x
=\int_{Y}u_{\rm per}\cdot\nabla\phi\,\dd x=\int_{Y}u\cdot\nabla\phi\,\dd x.
\end{align*}
\end{proof}
\begin{lmm}\label{efjohrioho1}
Let $u\in L^2(Y,\mathbb{R}^3)$. Assume that
for all $Y$-periodic $\phi\in\mathcal{C}^{\infty}(\mathbb{R}^3,
\mathbb{R})$
\begin{align*}
&\int_{Y}u\cdot\nabla\phi\,\dd x=0.
\end{align*}
Then for all $c=(c_1,c_2,c_3)\in\mathbb{R}^3$
and
$\psi\in\mathcal{C}^{\infty}_0(c+Y,\mathbb{R})$   
\begin{align*}
&\int_{c+Y}u_{\rm per}\cdot\nabla\psi\,\dd x=0.
\end{align*}
\end{lmm}
\begin{proof}
Extend $\psi$ periodically to $\mathbb{R}^3$ and denote the extension  by $\psi_{\rm per}$. 
By Lemma~\ref{elvhfbvhk1} we have
\begin{align*}
&\int_{c+Y}u_{\rm per}\cdot\nabla\psi\,\dd x=
\int_{c+Y}u_{\rm per}\cdot\nabla\psi_{\rm per}\,\dd x=
\int_{Y}u\cdot\nabla\psi_{\rm per}\,\dd x=0.
\end{align*}
\end{proof}
\begin{lmm}\label{thm135}
Let $u\in L^2(Y,\mathbb{R}^3)$.  Assume that
for all $Y$-periodic $\phi\in\mathcal{C}^{\infty}(\mathbb{R}^3,\mathbb{R})$
\begin{equation*}
\int_{Y}u\cdot\nabla\phi\,\dd x=0.
\end{equation*}
Then for all $\psi\in\mathcal{C}^{\infty}_0(\mathbb{R}^3,\mathbb{R})$
we have
\begin{align*}
&\int_{\mathbb{R}^3}u_{\rm per}\cdot\nabla\psi\,\dd x=0.
\end{align*}
\end{lmm}
\begin{proof}
Take any $\psi\in\mathcal{C}^{\infty}_0(\mathbb{R}^3,\mathbb{R})$.
Since $\supp\, \psi$ is compact, there exists $m\in \mathbb{N}$ such that
$\supp\, \psi\subset \prod_{i=1}^3 [-ma_i,ma_i]$.
For $k\in\mathbb{R}^3$, let us use the notation  $\vert k\vert=\max(\vert k_i\vert,\; i=1,2,3)$.
Let us denote for $k,j\in\mathbb{Z}^3$
\begin{align*}
Y_{k}&=\left(k_1a_1,k_2a_2,k_3a_3\right)+Y,\\
Y_{k,j}&=\left(\left(k_1+\frac{j_1}{2}\right)a_1,
\left(k_2+\frac{j_2}{2}\right)a_2,
\left(k_3+\frac{j_3}{2}\right)a_3\right)+Y.
\end{align*}
Note that 
$\supp\, \psi\subset \cup_{\vert k\vert\le m+1}\overline{Y}_k$ and 
$\overline{Y}_k\subset
\cup_{\vert j\vert\le 1}Y_{k,j}$.
Thus, the system $Y_{k,j}$, $\vert k\vert\le m+1$, 
$\vert j\vert\le 1$, 
is a finite open covering of $\supp\, \psi$.
Let us consider a partition of unity 
(see Theorem 5.3.8 in \cite{KFJ}) $\varphi_{k,j}$, 
$\vert k\vert\le m+1$, $\vert j\vert\le 1$, fulfilling:
$\varphi_{k,j}\in\mathcal{C}^\infty(\mathbb{R}^3,\mathbb{R})$,
$\supp\, \varphi_{k,j}\subset{Y}_{k,j}$, and 
$\sum_{|k|\le m+1,\vert j\vert\le 1}\varphi_{k,j}(x)=1$ for all $x\in \supp\, \psi$.
Set $\psi_{k,j}(x)=\psi(x)\varphi_{k,j}(x)$, $x\in\mathbb{R}^3$. 
Then $\psi_{k,j}\in \mathcal{C}^{\infty}_0(\mathbb{R}^3,\mathbb{R})
$, $\supp\, \psi_{k,j}\subset\overline{Y}_{k,j}\cap \supp\, \psi$, and
\begin{align*}
\psi(x)=\sum_{|k|\le m+1,\; \vert j\vert\le 1}\psi_{k,j}(x),
\quad x\in \supp\,\psi.
\end{align*}
Now we can write by Lemma \ref{efjohrioho1}
\begin{align*}
&\int_{\mathbb{R}^3}u_{\rm per}\cdot\nabla \psi\,\dd x=
\int_{\supp\, \psi}u_{\rm per}\cdot \sum_{|k|\le m+1,\; \vert j\vert\le 1}\nabla\psi_{k,j}\,\dd x
=\sum_{|k|\le m+1,\; \vert j\vert\le 1} 
\int_{\overline{Y}_{k,j}} u_{\rm per}\cdot\nabla {\psi}_{k,j}\,\dd x=0.
\end{align*}
\end{proof}

{\bf Proof of Theorem~\ref{charU}:}
\begin{proof}
Let us prove $ H^{{\rm div},0}_{\rm per}(Y,{\mathbb R}^3)\subset W$ first.
Let us have $u\in H^{{\rm div},0}_{\rm per}(Y,{\mathbb R}^3)$
and $\phi\in H^1_{\rm per}(Y,{\mathbb R})$. 
Then using Lemma~\ref{Pok1}
\begin{equation*}
\int_Y u\cdot\nabla\phi\,\dd x=\int_{\partial Y}(u\cdot n)\phi\, \dd s=0,
\end{equation*}
because the normal fluxes $u\cdot n$ have opposite signs on any pair of 
opposite sides of $Y$, while $\phi$ has the same traces there.
To prove $W\subset H^{{\rm div},0}_{\rm per}(Y,{\mathbb R}^3)$,
consider $u\in W$. Since
$C_0^\infty(Y,\mathbb{R})\subset H^1_{\rm per}(Y,{\mathbb R})$, 
we get 
\begin{equation*}
\int_Y u\cdot\nabla \phi\, \dd x=0
\end{equation*}
 for all $\phi\in C^\infty_{0}(Y,{\mathbb R})$, and thus
$u\in H^{{\rm div},0}(Y,{\mathbb R}^3)$. It remains to prove 
that $u$ is $Y$-periodic, i.e.~$u_{\rm per}\in H^{{\rm div},0}_{\rm loc}({\mathbb R}^3,{\mathbb R}^3)$.
For $u\in W$ 
and $\psi\in C_0^\infty(\mathbb{R}^3,\mathbb{R})$, we get from Lemma~\ref{thm135}
\begin{equation*}
\int_{\mathbb{R}^3}u_{\rm per}\cdot\nabla\psi\,\dd x=0.
\end{equation*}
\end{proof}

\section{Proof of Theorem~\ref{lem55}}\label{apxB}

\begin{lmm}\label{lem2}
Let $H$ be a real Hilbert space with the inner product
$(\cdot,\cdot)$ and let $H=V\oplus W$ where $V$ and $W$ are
nontrivial closed subspaces of $H$ and for all
$v\in V$ and $w\in W$ it holds $(v,w)=0$.
Let $F$ be a quadratic
functional on $H$ defined by
\begin{equation}\nonumber
F(u)=(A(u-u_0),u-u_0)
\end{equation}
where $A:H\to H$ is a symmetric positive definite operator and $u_0\in H$ and fixed.
Then 
\begin{equation}\nonumber
\min_{v\in V}F(v)=\max_{w\in W}\left(
\min_{u\in H}\left( F(u)-(w,u)\right)\right).
\end{equation}
\end{lmm}
\begin{proof}
We prove "$\ge$" first. We have for all $w\in W$
\begin{eqnarray}
\min_{v\in V}\left( F(v)-(w,v)\right)&\ge& \min_{u\in H}\left( F(u)-(w,u)\right)\nonumber\\
\min_{v\in V} F(v)&\ge& \min_{u\in H}\left( F(u)-(w,u)\right)\nonumber
\end{eqnarray}
where we use $(v,w)=0$ for $v\in V$ and $w\in W$. Thus
\begin{equation}\nonumber
\min_{v\in V}F(v)\ge \max_{w\in W}\left(
\min_{u\in H}\left( F(u)-(w,u)\right)\right).
\end{equation}
To prove "$\le$" let us first note that the minimum 
of $F(v)$ in $V$ is achieved for such $v_0\in V$ that
\begin{equation}\label{v1}\tag{b.1}
(Av_0-Au_0,v)=0,\quad\text{for all}\;v\in V.
\end{equation}
Then using~\eqref{v1} for $v=v_0$ we get
\begin{equation}\label{v2}\tag{b.2}
\min_{v\in V}F(v)=F(v_0)=(A(v_0-u_0),v_0-u_0)=(Au_0,u_0)-(Av_0,v_0).
\end{equation}
Let us take $w_0=2A(v_0-u_0)$. From~\eqref{v1} we have $w_0\in W$.
We have
\begin{eqnarray}
F(u)-(w_0,u)&=&(A(u-u_0),u-u_0)-2(Av_0-Au_0,u)\nonumber\\
&=&(A(u-v_0),u-v_0)+(Au_0,u_0)-(Av_0,v_0).\nonumber
\end{eqnarray}
Clearly, the minimum of $F(u)-(w_0,u)$
in $H$ is obtained for $u=v_0$, and then
\begin{equation}\nonumber
\min_{u\in H}\left( F(u)-(w_0,u)\right)=(Au_0,u_0)-(Av_0,v_0),
\end{equation}
which is equal to~\eqref{v2}. This proves "$\le$".
\end{proof}
\begin{lmm}\label{Lem1}
Let $\xi\in{\mathbb R}^3$ be arbitrary and fixed. Let
\begin{equation}\nonumber
(  A^*\xi,\xi)_{{\mathbb R}^3}=\frac{1}{\vert Y\vert}\min_{u\in U}
\int_Y(  A(\xi+\nabla u),\xi+\nabla u)_{{\mathbb R}^3}\, \dd x.
\end{equation}
Then we have
\begin{equation}\nonumber
(  A^*\xi,\xi)_{{\mathbb R}^3}=\frac{1}{\vert Y\vert}
\max_{w\in  W}
\int_Y 
-(A^{-1}w,w)_{\mathbb{R}^3}+
2(w,\xi)_{{\mathbb R}^3} \, \dd x.
\end{equation}
\end{lmm}
\begin{proof}
Let $\xi\in{\mathbb R}^3$ be fixed.
Let us define
\begin{equation}\nonumber
F(v)=\frac{1}{2}\int_Y(  A(\xi+v),
\xi+v)_{{\mathbb R}^3}\, \dd x,\quad v\in V.
\end{equation}
Then from Lemma~\ref{lem2} we get 
\begin{eqnarray}\nonumber
\frac{\vert Y\vert}{2}(  A^*\xi,\xi)_{{\mathbb R}^3}&=&\min_{v\in{V}}F(v)=
\max_{w\in W}\left(\min_{u\in{H}}\left(F(u)-
(w,u)_{H}\right)\right)\nonumber\\
&=& \max_{w\in W}\left(\min_{u\in{H}}\left(
\int_Y\frac{1}{2}(  A(\xi+u),
\xi+u)_{{\mathbb R}^3}\, \dd x-
\int_Y(  w,u)_{{\mathbb R}^3}\, \dd x\right)\right)\nonumber\\
&=& \max_{w\in W}\left(\min_{u\in{H}}
\int_Y\frac{1}{2}(  A(\xi+u),
\xi+u)_{{\mathbb R}^3}-(  w,u)_{{\mathbb R}^3}\, \dd x\right).\nonumber
\end{eqnarray} 
For any $w\in W$ the minimum is attained when $u=-\xi+A^{-1}w$
a.e.~in $Y$. Then
\begin{eqnarray}\nonumber
\frac{\vert Y\vert}{2}(  A^*\xi,\xi)_{{\mathbb R}^3}&=&
\max_{w\in W}\left(
\int_Y-\frac{1}{2}(  A^{-1}w,
w)_{{\mathbb R}^3} +(  w,\xi )_{{\mathbb R}^3}\, \dd x\right).\nonumber
\end{eqnarray} 
\end{proof}

{\bf Proof od Theorem~\ref{lem55}:}
\begin{proof}
Let $\xi\in\mathbb{R}^3$ be arbitrary and let
$w\in W$ be split into $w=\beta+v$, 
$v\in W_0$, $\beta\in\mathbb{R}^3$.
From Lemma~\ref{Lem1} and the definition~\eqref{min2} we have 
\begin{eqnarray}
(A^*\xi,\xi)_{{\mathbb R}^3}&=&\frac{1}{\vert Y\vert}\max_{w\in W}
2\int_Y (w,\xi)_{{\mathbb R}^3}\, \dd x-\int_Y(A^{-1}w,w)_{{\mathbb R}^3}\, \dd x\nonumber\\
&=&\frac{1}{\vert Y\vert}\max_{v\in W_0,\beta\in {\mathbb R}^3}
2\int_Y (\beta+v,\xi)_{{\mathbb R}^3}\, \dd x-\int_Y(A^{-1}(\beta+v),\beta+v)_{{\mathbb R}^3}\, \dd x\nonumber\\
&=&\frac{1}{\vert Y\vert}\max_{v\in W_0,\beta\in {\mathbb R}^3}
2\int_Y (\beta,\xi)_{{\mathbb R}^3}\, \dd x-\int_Y(A^{-1}(\beta+v),\beta+v)_{{\mathbb R}^3}\, \dd x\nonumber\\
&=&\frac{1}{\vert Y\vert}\max_{\beta\in {\mathbb R}^3}
\left(2\int_Y (\beta,\xi)_{{\mathbb R}^3}\, \dd x-
\min_{v\in W_0}\int_Y(A^{-1}(\beta+v),\beta+v)_{{\mathbb R}^3}\, \dd x\right)\nonumber\\
&=&\max_{\beta\in {\mathbb R}^3}
\left(\frac{2}{\vert Y\vert}\int_Y (\beta,\xi)_{{\mathbb R}^3}\, \dd x-
(B^*\beta,\beta)_{{\mathbb R}^3}\right)\nonumber\\
&=&\max_{\beta\in {\mathbb R}^3}
\left(2(\beta,\xi)_{{\mathbb R}^3}-(B^*\beta,\beta)_{{\mathbb R}^3}\right)=((B^*)^{-1}\xi,\xi)_{{\mathbb R}^3}.\nonumber
\end{eqnarray}
Thus we have $(B^*)^{-1}=A^*$.
Let us assume that the minimum in~\eqref{min1} is
achieved for $u=u_\xi\in U$.
We prove that $w_1$ defined by~\eqref{AEj}
is in $W$. Indeed, $w_1=A\xi+A\nabla u_\xi-A^*\xi$ is in $H$ and for any $\phi\in U$
\begin{equation}
\int_Y w_1\cdot\nabla \phi\, \dd x=
\int_Y  (A\xi+A\nabla u_\xi-A^*\xi)\cdot\nabla \phi\, \dd x
=\int_Y  (A\xi+A\nabla u_\xi)\cdot\nabla \phi\, \dd x=0,
\nonumber
\end{equation}
where the last equality follows from~\eqref{form111}. Thus $w_1\in W$.
By~\eqref{min1} and~\eqref{min11} for any $\mu\in\mathbb{R}^3$ we obtain 
\begin{eqnarray}
\int_Y(w_1,\mu)_{\mathbb{R}^3}\, \dd x&=&\int_Y( A(\xi+\nabla u_\xi),\mu)_{\mathbb{R}^3}\,\dd x
-\vert Y\vert\, (A^*\xi,\mu)_{\mathbb{R}^3}\nonumber\\
&=&
\int_Y( A(\xi+\nabla u_\xi),\mu)_{\mathbb{R}^3}\,\dd x
-\int_Y( A(\xi+\nabla u_\xi),\mu)_{\mathbb{R}^3}\,\dd x=0.\nonumber
\end{eqnarray}
Thus $w_1\in W_0$.
Recalling that $\alpha=A^*\xi$ and the minimizing property of $u=u_\xi$ in~\eqref{min1} 
we have
\begin{eqnarray}
(B^*\alpha,\alpha)_{{\mathbb R}^3}&=&(A^*\xi,\xi)_{{\mathbb R}^3}=\frac{1}{\vert Y\vert}\int_Y({A}\,(\xi+\nabla u_\xi),\xi+\nabla u_\xi)_{{\mathbb R}^3}
\, \dd x\nonumber\\
&=&\frac{1}{\vert Y\vert}\int_Y({A^{-1}}\,(A\xi+A\nabla u_\xi),A\xi+A\nabla u_\xi)_{{\mathbb R}^3}
\, \dd x\nonumber\\
&=&\frac{1}{\vert Y\vert}\int_Y({A^{-1}}\,(\alpha+w_1),\alpha+w_1)_{{\mathbb R}^3}
\, \dd x.\nonumber
\end{eqnarray}
Since the minimizing field $w_\alpha$ is unique,
we get $w_1=w_\alpha$.
\end{proof}

\end{document}